\documentclass[preprint,12pt]{elsarticle}

\usepackage{graphicx}
\usepackage{wrapfig}
\usepackage{subfigure}
\usepackage{multicol}
\usepackage{natbib}
\usepackage{multirow}
\usepackage{stfloats}

\newcommand{\tr}{\operatorname{tr}}
\newcommand{\nn}{\nonumber}

\newcommand{\rank}{\mathrm{rank}}
\newcommand{\Es}{\mathbb{E}}
\usepackage{CJK}
\usepackage{amsmath}
\usepackage{bbm}
\usepackage{amsfonts}
\usepackage{amssymb}
\usepackage{enumerate}
\usepackage{epstopdf}
\usepackage{algpseudocode}
\usepackage{algorithm}
\usepackage{amsfonts}
\usepackage[dvipsnames]{xcolor}
\newcommand{\mz}{\color{black}}
\newcommand{\sy}{\color{black}}

\newtheorem{thm}{Theorem}
\newtheorem{lem}[thm]{Lemma}
\newtheorem{prob}{Problem}
\newtheorem{corollary}{Corollary}
\newtheorem{proposition}{Proposition}
\newdefinition{rmk}{Remark}
\newproof{pf}{Proof}
\newproof{pot}{Proof of Theorem \ref{thm2}}

\date{}
\begin{document}

\begin{frontmatter}

\title{Robust fixed-lag smoothing under model perturbations}

{\mz \author[Italy]{Shenglun Yi}\ead{yishenglun@dei.unipd.it}}   
\author[Italy]{Mattia Zorzi}\ead{zorzimat@dei.unipd.it}               

\address[Italy]{Department of Information Engineering, University of Padova, Via Gradenigo 6/B, 35131 Padova, Italy}             


\begin{abstract}
A robust fixed-lag smoothing approach is proposed in the case there is a  mismatch between the nominal model and the actual model. The resulting robust smoother is characterized by a dynamic game between two players: one player selects the least favorable model in a prescribed ambiguity set, while the other player selects the fixed-lag smoother minimizing the smoothing error with respect to least favorable model. We propose an efficient implementation of the proposed smoother. Moreover, we characterize the corresponding least favorable model over a finite time horizon. Finally, we test the robust fixed-lag smoother in two examples. The first one regards a target tracking problem, while the second one regards a parameter estimation problem.
\end{abstract}

\end{frontmatter}


\section{Introduction}

Fixed-lag smoothing aims to estimate the state of a dynamical system at time $t$ using the observations in the interval $[0,t+L-1]$ with $L>1$. This algorithm is fundamental in various applications, e.g. tracking and navigation because it can handle online requirements, see \cite{cattivelli2010diffusion,dong2011motion,papi2014fixed}. Moreover,
it can be used in the expectation maximization (EM) algorithm  to compute  the maximum likelihood (ML) estimator of the unknown parameters characterizing the matrices of a state space model, see e.g. \cite{shumway1982approach,5717145,sarkka2013bayesian}.
Indeed, although the EM algorithm is based on the Rauch-Tung-Striebel (RTS) smoother, \cite{linearestimationM}, its estimate can be approximated by the one given by the fixed-lag smoother provided that $L$ is taken sufficiently large. However, in all the aforementioned applications, the actual model is typically known only imprecisely that is only the nominal model is known. In this situation, these smoothers could perform poorly.

Model uncertainty is traditionally addressed by risk sensitive filtering, see e.g. \cite{RISK_WHITTLE_1980,RISK_PROP_BANAVAR_SPEIER_1998,H_INF_HASSIBI_SAYED_KAILATH_1999,huang2018distributed,LEVY_ZORZI_RISK_CONTRACTION,OPTIMAL_SPEYER_FAN_BANAVAR_1992} . Here, the robust estimator minimizes an exponential loss function which severely penalizes large errors. These filters can be also interpreted as the solution of a dynamic minimax game  \cite{yoon2004robust}, see also {\mz \cite{boel2002robustness,HANSEN_SARGENT_2005, HANSEN_SARGENT_2007, zorzi2018robust}}. Then, \cite{ROBUST_STATE_SPACE_LEVY_NIKOUKHAH_2013} proposed a robust estimator which solves an incremental minimax game. At time $t$ the actual model belongs to the ambiguity set which is a ball, in the Kullback-Leibler (KL) topology, about the nominal model. In this way, the uncertainty is ``spread'' along the time and not concentrated in specific time steps. Then, there are two players which operate against. One player, say nature, selects the least favorable model in this prescribed ``ball", and the other player designs the optimum estimator for the least favorable model. It is worth noting that many extensions of this paradigm have been proposed such as: the case with different ambiguity sets  {\mz  \cite{abadeh2018wasserstein,STATETAU_2017,zorzi2017convergence,zorzi2017robustness}}; the distributed case \cite{zorzi2019distributed,RKDISTR_opt}; the case with external input \cite{RS_MPC_IET}; the case of degenerate densities \cite{yi2020lowTAC,yi2020low}.

In the literature, the robust smoothing problems mainly consider two situations. In the first case, the noise distribution is known but it is not {\sy necessarily} Gaussian {\sy \cite{piche2012recursive,aravkin2011ell,farahmand2011doubly,huang2016robust,navarro2019widely, fernandez2020semi, zhao2021discrete, huang2019robust,zhao2017bayesian,huang2020slide}}, for instance the noise process is assumed to have a non-Gaussian distribution in order to model outliers, {\sy temporary model uncertainties,} missing observations or sensor delays. {\sy Some of these robust paradigms are adaptive because the parameters of the noises characterizing  the state space model are inferred from the collected data.}
In the second situation, the noise distribution is not known but this process takes values in a bounded set, e.g. an ellipsoidal set \cite{kwon2015robust,rehman2016robust}. However, there are relatively few studies on robust smoothing problems which use the risk sensitive philosophy, see \cite{theodor1994game,bolzern2004discrete}.

In this paper, we propose a new robust fixed-lag smoothing problem where the model uncertainty is expressed incrementally as in \cite{ROBUST_STATE_SPACE_LEVY_NIKOUKHAH_2013,yi2020lowTAC}. Thus, at each time step we have to solve a dynamic game between two players: the nature which selects the least favorable model in the ambiguity set and the other player which designs the optimal fixed-lag smoother according to the least favorable model. The resulting smoother is characterized  by matrices whose dimension is proportional to the lag $L$. On the other hand, the typical value of the lag is large. Accordingly, numerical instabilities and high computational burden characterize the algorithm. In order to avoid those issues, we propose an efficient implementation drawing inspiration from the reduced order fixed-lag smoother proposed in \cite{discretetimefixedlag}. Then, the corresponding least favorable model over a finite simulation horizon is derived to evaluate the performance of the smoother. Finally, we consider a target tracking problem and a parameter estimation problem to test the performance of the robust fixed-lag smoother.

The outline of the paper is as follows. In Section \ref{sec_2}, we introduce the problem formulation about robust fixed-lag smoothing. In  Section  \ref{sec_3} we derive the robust fixed-lag smoother. The corresponding algorithm is then reformulated in an efficient way to reduce the computational burden in Section \ref{sec_4}. In Section \ref{sec_5}, we derive the least favorable model corresponding to the robust fixed-lag smoother.  The numerical examples are provided in Section \ref{sec_7}, which is devoted to the target tracking problem, and in Section \ref{sec_6},  which is devoted to the parameter estimation problem.
 Finally, in Section \ref{sec_8} we draw the conclusions.

\section{Problem Formulation} \label{sec_2}

We consider the nominal state space model:
\begin{equation}\label{ss_space}
\begin{array}{cc}x_{t+1}&=A x_t+B v_t \\ y_t&=C x_t+D v_t\end{array}
\end{equation}
where {$A \in \mathbb{R}^{n\times n}$, {$B \in \mathbb{R}^{n\times (m+n)}$}, $C \in \mathbb{R}^{m\times n}$ and $D \in \mathbb{R}^{m\times (m+n)}$, $x_t$ is the state vector, $y_t $ is the observation vector, and {$v_t \in \mathbb{R}^{m+n}$} is normalized white Gaussian noise. Moreover, $x_0 \sim \mathcal{N}\left(\hat x_0, V_0\right)$, with $V_0>0$},  which is independent from $v_t$. We also assume that $BD^\top=0$, $\mathrm{rank}(B)=n$ and $\mathrm{rank}(D)=m$. {\sy In this way, the process noise, {\sy say $\mathrm w_t:=B v_t$,} and the measurement noise, {\sy say $\mathrm v_t:=D v_t$, are independent and their covariance matrices are  $BB^{\top}$ and $DD^{\top}$}, respectively}. Fixed-lag smoothing aims to find an estimate of {$\hat x_{t-L+1|t}$} of {$x_{t-L+1}$} given $Y_t=\left\{{y_0 }\cdots y_t \right\}$ and $L$ denotes the lag. It is well-known that such a problem can be interpreted as a Kalman prediction problem corresponding to the following nominal augmented state space model \cite{discretetimefixedlag}:
\begin{equation}
\begin{array}{rl}\xi_{t+1}&=\tilde{A} \xi_t+\tilde{B} v_t \\ y_t&=\tilde{C} \xi_t+\tilde{D} v_t  \\  x_{t-L+1}&=\tilde H\xi_{t+1}\end{array}  \label{model}
\end{equation}
where $\tilde A \in \mathbb{R}^{(L+1)n\times (L+1)n}$, $\tilde B \in \mathbb{R}^{(L+1)n\times (m+n)}$, $\tilde C \in \mathbb{R}^{m\times (L+1)n}$, $\tilde D \in \mathbb{R}^{m\times (m+n)}$ and $\tilde H \in \mathbb{R}^{n\times(L+1)n}$ are such that
$$\xi_{t}=\left[\begin{array}{c}{x_{t}} \\ {x_{t-1}} \\ {\vdots} \\ {x_{t-L}}\end{array}\right],~
\tilde{A}=\left[\begin{array}{cccc}{A} & {0} & {\cdots} & {0} \\ {I} & {\cdots} & {0} & {0} \\ {\vdots} & {\ddots} & {\vdots} & {\vdots} \\ {0} & {\cdots} & {I} & {0}\end{array}\right],~ \tilde{B}=\left[\begin{array}{c}{B} \\ {0} \\ {\vdots} \\ {0}\end{array}\right]$$ $$\tilde{C}=\left[\begin{array}{llll}{C} & {0} & {\cdots} & {0}\end{array}\right], \tilde{D}=D, \tilde H=\left[\begin{array}{llll}{0} & {0} & {\cdots} & {I}\end{array}\right].$$
Then, $${ \hat x_{t-L+1|t}}=\tilde H\hat \xi_{t+1}$$ where $\hat \xi_{t+1}$ is the one step-ahead predictor of $\xi_{t+1}$ given $Y_t$.  We define $z_t:=[\,\xi_{t+1}^{\top } \; y_t^{\top }\,]^{\top }.$  Let $\phi (z_t|\xi_t)$ be the transition probability density function of $z_t$ given $\xi_t$ corresponding to the nominal model (\ref{model}). Then, $\phi_t(z_t|\xi_t) \sim \mathcal{N}\left(m_{z_t}, K_{z_t}\right)$ with$$
m_{z_t}=\left[\begin{array}{c}{\tilde A} \\ {\tilde C}\end{array}\right] \xi_{t}, ~~~K_{z_t}=\left[\begin{array}{cc}{\tilde B \tilde B^\top} & {0} \\ {0} & {\tilde D \tilde D^\top}\end{array}\right].$$
Notice that $\phi_t(z_t|\xi_t)$ is a degenerate probability density function because $\tilde B \tilde B^\top$ is singular and thus $K_{z_t}$ as well. More precisely, $\rank(K_{z_t})=n+m$. Accordingly, the support of $\phi_t(z_t|\xi_t)$ is the {$n+m$}-dimensional affine subspace
\begin{equation*}
\mathcal{A}_t=\left\{m_{z_t}+v, \quad v \in \mathrm {Im} \left(K_{z_t}\right)\right\}
\end{equation*} which depends on $\xi_t$.
Then,
\begin{equation*}
\begin{aligned} \phi_t(z_t|\xi_t) =\left[(2 \pi)^{n+m} \operatorname{det}^{+}\left(K_{z_t}\right)\right]^{-1 / 2}  \times  \exp \left[-\frac{1}{2}\left(z_t-m_{z_t}\right)^{\top } K_{z_t}^{+}\left(z_t-m_{z_t}\right)\right] \end{aligned}
\end{equation*}
where $K^{+}_{z_t}$ is the pseudo-inverse of $K_{z_t}$ and $\operatorname{det}^{+}(K_{z_t})$ is the pseudo-determinant of $K_{z_t}$.

The nominal model in (\ref{model}) in the time horizon $[0,N]$ is described by the joint probability density
\begin{equation}
f\left(\Xi_{N+1}, Y_{N}\right)={\tilde f_{0}\left(\xi_{0}\right)} \prod_{t=0}^{N} \phi_{t}\left(z_{t} | \xi_{t}\right) \label{f(T)}
\end{equation}
where
\begin{align*}
\Xi_{N+1}^\top&=\left[\begin{array}{ccccc}{\xi_{0}^\top} & {\ldots} & {\xi_{t}^\top} & {\ldots} & {\xi_{N+1}^\top}\end{array}\right]\\ Y_{N}^\top&=\left[\begin{array}{ccccc}{y_{0}^\top} & {\ldots} & {y_{t}^\top} & {\ldots} & {y_{N}^\top}\end{array}\right],
\end{align*}
${\tilde f_{0} (\xi_{0})} \sim \mathcal{N} (\hat \xi_{0}, \tilde V_0 )$ with
\begin{align}\label{set_init} \hat \xi_0=\left[\begin{array}{c}\hat x_0 \\  \star \end{array}\right], ~\tilde V_0=\left[\begin{array}{cc}
V_0 & 0 \\ 0 & \star \end{array}\right].\end{align} In the above equations the star symbol means that it is an arbitrary vector or matrix. Indeed, since $x_0$ is independent from $x_t$ with $-L\leq t<0$, the smoother does not depend on those parameters for $t\geq 0$. Accordingly, without loss of generality, we assume that $\tilde V_0>0$.

In this paper, we consider the situation in which the actual model does not coincide with the nominal one in (\ref{model}). In particular, we assume that the probability density of the actual model has a structure similar to the one in (\ref{f(T)}):
\begin{equation*}
\tilde{f}\left(\Xi_{N+1}, Y_{N}\right)=\tilde{f}_{0}\left(\xi_{0}\right) \prod_{t=0}^{N} \tilde{\phi}_{t}\left(z_{t} | \xi_{t}\right)
\end{equation*}
where we assume that $\tilde \phi_t$ has the same support of $\phi_t$. In this way, we can measure the discrepancy between $f$ and $\tilde f$ through the KL-divergence:
\begin{equation*}
\begin{aligned} &D(\tilde f_Z, f_Z) := \int_{\mathcal A_\Xi}\int_{\mathcal A_Y}\tilde{f}\left(\Xi_{N+1}, Y_{N}\right) \ln \frac{\tilde{f}\left(\Xi_{N+1}, Y_{N}\right)}{f\left(\Xi_{N+1}, Y_{N}\right)}   d Y_N d \Xi_{N+1} \end{aligned}
\end{equation*}
{\sy where $\tilde E[\cdot]$ is the expected value operator with respect to the actual probability density $\tilde f$ and $\mathcal A_\Xi\times \mathcal A_Y$ is the support of $f$ and $\tilde f$. It is worth noting that the KL divergence is the natural metric to measure such a mismatch in the case the nominal model is inferred from data, see \cite{ROBUSTNESS_HANSENSARGENT_2008}.}
It is not difficult to see that
\begin{equation}
\begin{aligned}
 D(\tilde{f}, f)= \sum_{t=0}^{N} D\left(\tilde{\phi}_{t}, \phi_{t}\right)\label{D}
\end{aligned}
\end{equation}
where
\begin{equation*}
\begin{aligned} D (&\tilde{\phi}_{t}, \phi_{t} ) =\tilde{\Es}\left[\ln \left(\frac{\tilde{\phi}_{t}\left(z_{t} | \xi_{t}\right)}{\phi_{t}\left(z_{t} | \xi_{t}\right)}\right)\right] \\ &:= \int_{\bar{\mathcal A}_t }\int_{\mathcal A_t} \tilde{\phi}_{t}\left(z_{t} | \xi_{t}\right) \tilde{f}_{t}\left(\xi_{t}\right) \ln \left(\frac{\tilde{\phi}_{t}\left(z_{t} | \xi_{t}\right)}{\phi_{t}\left(z_{t} | \xi_{t}\right)}\right) d z_{t} d \xi_{t}\end{aligned}
\end{equation*}
where $\bar{\mathcal A}_t$ is the support of $\tilde f_t(\xi_t)$  which denotes the actual marginal density of $\xi_t$.

Since the actual model is not known, we assume the latter belongs to the ambiguity set which is a ball about $f$ formed by placing an upper bound on $D(\tilde{f}, f)$. However, this ambiguity set contains models which concentrate the uncertainty in a unique time step, i.e. a situation which is unrealistic in practice. On the other hand, in view of Equation (\ref{D}), we can express such a mismatch incrementally through $\phi_{t}$ and $\tilde{\phi}_{t}$. Accordingly, we assume that $\tilde \phi_t$ given $Y_{t-1}$ belongs to the following ambiguity set:
\begin{equation*}
\mathcal B_t:=\left\{\, \tilde{\phi}_{t} \hbox{ s.t. }\tilde{\Es}\left[\ln \left(\frac{\tilde{\phi}_{t}\left(z_{t} | \xi_{t}\right)}{\phi_{t}\left(z_{t} | \xi_{t}\right)}\right) \bigg| Y_{t-1}\right] \leq c_{t}\right\}
\end{equation*}
where
\begin{equation}
\begin{aligned}&\tilde{\Es}\left[\ln \left(\frac{\tilde{\phi}_{t}\left(z_{t} | \xi_{t}\right)}{\phi_{t}\left(z_{t} | \xi_{t}\right)}\right) \bigg| Y_{t-1}\right] \\ &:= \int_{ \check{\mathcal A}_t}\int_{\mathcal A_t} \tilde{\phi}_{t}\left(z_{t} | \xi_{t}\right) \tilde{f}_{t}\left(\xi_{t} | Y_{t-1}\right) \ln \left(\frac{\tilde{\phi}_{t}\left(z_{t} | \xi_{t}\right)}{\phi_{t}\left(z_{t} | \xi_{t}\right)}\right) d z_{t} d \xi_{t},
\label{tildeE}
\end{aligned}
\end{equation}
and $\check A_t$ is the support of $\tilde f_t(\xi_t|Y_{t-1})$. It is worth noting that $c_t>0$, hereafter called tolerance, is the mismodeling budget allowed at time step $t$. Our aim is to address the following problem.

\begin{prob}
Design a fixed-lag smoother with respect to  the ambiguity set $\mathcal B_t$ for $t=0\ldots N$.
\end{prob}

\section{ Robust smoothing}\label{sec_3}

We propose a robust fixed-lag smoother of $x_{t-L+1}$ given $Y_t$  with respect to $\mathcal B_t$ solving the following minimax problem:
\begin{align}  &\hat x_{t-L+1|t}=\tilde H \hat \xi_{ t+1}\nn \\
& \label{minimax} \hat \xi_t =\underset{g_t \in \mathcal{G}_{t}}{\mathrm{argmin}}\max_{\tilde{\phi}_{t} \in \mathcal{B}_{t}} J_t(\tilde {\phi}_t,g_t)
\end{align}
where \begin{equation*}
\begin{aligned}
J_t(\tilde {\phi}_t,g_t)=&\frac{1}{2}\tilde{\mathbb{E}}\left[\|\tilde H\left(\xi_{t+1}-g_t\left(y_{t}\right)\right)\|^{2} | Y_{t-1}\right]\\
=&\frac{1}{2} \int_{\check A_t}\int_{\mathcal A_t}\|\tilde H\left(\xi_{t+1}-g_{t}\left(y_{t}\right)\right)\|^{2} \tilde{\phi}_{t}\left(z_{t} | \xi_{t}\right)  \times \tilde{f}_{t}\left(\xi_{t} | Y_{t-1}\right) d z_{t} d \xi_{t},
\end{aligned}
\end{equation*}
$\mathcal{G}_{t}$ denotes the class of estimators with finite second-order moments with respect to all the densities $\tilde{\phi}_{t}\left(z_{t} | \xi_{t}\right) \tilde{f}_t(\xi_t|Y_{t-1})$ such that $\tilde{\phi}_{t} \in \mathcal{B}_{t}$.  {Notice that $\tilde \phi_t$ must satisfy the constraint:}
\begin{equation}
I_{t}(\tilde{\phi}_{t}) \triangleq \int_{ \check A_t}\int_{\mathcal A_t} \tilde{\phi}_{t}\left(z_{t} | \xi_{t}\right) \tilde{f}_{t}\left(\xi_{t} | Y_{t-1}\right) d z_{t} d \xi_{t}=1. \label{I}
\end{equation}
{\mz It is worth noting that problems like (\ref{minimax}) can be written as a risk-sensitive problem, i.e. as a minimization problem where the standard quadratic cost function is replaced by an exponential cost function, see \cite{levy2004robust} for more details.}

\begin{lem}
\label{lemma1}
 For a fixed estimator $g_t \in \mathcal{G}_t$, the density $\tilde {\phi}_{t}\left(z_{t} | \xi_{t}\right) \in \mathcal{B}_{t}$ that maximizes the objective function
$$J_t(\tilde {\phi}_t,g_t)=\tilde{\mathbb{E}}\left[\|\tilde H\left(\xi_{t+1}-g_t\left(y_{t}\right)\right)\|^{2} | Y_{t-1}\right]$$
under constraint $D_{t}(\tilde{\phi}_{t}, \phi_{t}) \leq c_{t}$ is given by
\begin{equation} \label{phi_0}
\tilde{\phi}_{t}^{0}=\frac{1}{M_{t}\left(\lambda_{t}\right)} \exp \left(\frac{1}{2 \lambda_{t}}\left\|\tilde H(x_{t+1}-g_{t}\left(y_{t})\right)\right\|^{2}\right) \phi_{t}
\end{equation}  where $M_t(\lambda_t) $ is the normalizing constant defined as follows:
\begin{equation*}
\begin{aligned} M_{t}\left(\lambda_{t}\right)=&\int_{\check{\mathcal A}_t}\int_{\mathcal A_t} \exp \left(\frac{1}{2 \lambda_{t}}\|\tilde H(\xi_{t+1}-g_{t}\left(y_{t})\right)\|^{2}\right) \phi_{t}  \times \tilde{f}_{t}\left(\xi_{t} | Y_{t-1}\right) d z_{t} d \xi_{t}. \end{aligned}
\end{equation*} Moreover, for $c_t>0$ sufficiently small, there exists a unique $\lambda_t>0$ such that $D(\tilde \phi^0_t, \phi_t)=c_t $.
\end{lem}

\smallskip
\begin{pf}  The proof is similar to the one of \cite[Lemma2]{yi2020lowTAC}. \qed\end{pf}

Once we get the function $\tilde \phi^0_t$, the estimator $g_t \in \mathcal G_t$ minimizing the objective function $J_t(\tilde \phi_t^0,g_t)$ is given by
\begin{equation*}
\begin{aligned} \hat{\xi}_{t+1}=g_{t}^{0}\left(y_{t}\right) =\tilde{\Es}\left[\xi_{t+1} | Y_{t}\right] =\int_{ \check  {\mathcal A}_{t+1}} \xi_{t+1} \tilde{f}_{t+1}\left(\xi_{t+1} | Y_{t}\right) d \xi_{t+1} \end{aligned}
\end{equation*}
where \begin{equation}
\tilde{f}_{t+1}\left(\xi_{t+1} | Y_{t}\right)=\frac{\int_{\check{\mathcal A_t}} \tilde{\phi}_{t}^{0}\left(z_{t} | \xi_{t}\right)  \tilde{f}_{t}\left(\xi_{t} | Y_{t-1}\right) d \xi_{t}}{\int_{\check{\mathcal A}_t} \int_{\mathcal A_{t}^\star} \tilde{\phi}_{t}^{0}\left(z_{t} | \xi_{t}\right)  \tilde{f}_{t}\left(\xi_{t} | Y_{t-1}\right) d \xi_{t+1} d \xi_{t}}
\label{tildef}
\end{equation}
where $\mathcal A^\star_t$ is defined as follows: $\xi_{t+1}\in \mathcal A^\star_t$ if and only if there exists at least one $y_t$ for which $ [\,\xi_{t+1}^\top\; y_t^\top\,]^\top\in\mathcal A_t $.
The optimal estimator $g_t^0$ {solution to (\ref{minimax})} relies on the least-favorable density $\tilde{\phi}_{t}^{0}\left(z_{t} | \xi_{t}\right)$. On the other hand, the latter depends on the estimator $g^0_t$. In order to break this deadlock problem, an additional assumption is needed. More precisely, we assume that the $a ~priori$ conditional density $\tilde{f}_{t}(\xi_{t} | Y_{t-1})$ is Gaussian  $\tilde{f}_{t}(\xi_{t} | Y_{t-1}) \sim \mathcal{N}(\hat{\xi}_{t}, \tilde V_{t})$. In view of (\ref{model}), the marginal density
\begin{equation*}
\bar{f}_{t}\left(z_{t} | Y_{t-1}\right):=\int_{\mathcal A_t} \phi_{t}\left(z_{t} | \xi_{t}\right) \tilde{f}_{t}\left(\xi_{t} | Y_{t-1}\right) d \xi_{t}
\end{equation*}
is Gaussian so that
\begin{equation}
\bar{f}_{t}\left(z_{t} | Y_{t-1}\right) \sim \mathcal{N}\left(m_{z_{t} | Y_{t-1}}, K_{z_{t} | Y_{t-1}}\right) \label{barf}
\end{equation}
where $$\begin{aligned} m_{z_{t} | Y_{t-1}}& = \left[\begin{array}{c}{\tilde{A}} \\ {\tilde{C}}\end{array}\right] \hat{\xi}_{t},\\
K_{z_{t} | Y_{t-1}}&=\left[\begin{array}{c}{\tilde{A}} \\ {\tilde{C}}\end{array}\right]\tilde V_{t}\left[\begin{array}{cc}{\tilde{A}^{\top }} & {\tilde{C}^{\top }}\end{array}\right]+\left[\begin{array}{c}{\tilde{B}} \\ {\tilde{D}}\end{array}\right]\left[\begin{array}{cc}{\tilde{B}^{\top }} & {\tilde{D}^{\top }}\end{array}\right].
\end{aligned}$$Then, on the basis of Lemma \ref{lemma1}, we have that the least favorable density of $z_t$ given $Y_{t-1}$ is
\begin{equation}
\begin{aligned}
&\tilde{f}_{t}\left(z_{t} | Y_{t-1}\right):=\int_{\mathcal A_t} \tilde{\phi}_{t}^{0}\left(z_{t} | \xi_{t}\right) \tilde{f}_{t}\left(\xi_{t} | Y_{t-1}\right) d \xi_{t}\\
&=\frac{1}{M_t\left(\lambda_{t}\right)} \exp \left(\frac{1}{2 \lambda_{t}}\|\tilde H\left(\xi_{t+1}-g_{t}\left(y_{t}\right)\right)\|^{2}\right)  \bar{f}_{t}\left(z_{t} | Y_{t-1}\right).
\label{tildef0}
\end{aligned}
\end{equation}
Accordingly, $\tilde{f}_{t}\left(z_{t} | Y_{t-1}\right)$ is a Gaussian probability density.

\begin{lem}\label{lemma3}
Consider the state space model (\ref{ss_space}) with $\mathrm{rank} (B)=n$ and $\mathrm{rank} (D)=m$. If $\tilde V_t>0$ then $K_{z_t|Y_{t-1}}>0$, i.e. $\bar{f}_{t}\left(z_{t} | Y_{t-1}\right)$ is a non-degenerate density. Moreover,
\begin{align}
 \tilde P_{t+1}:=\nonumber \tilde A\tilde V_t \tilde A^\top-\tilde A\tilde V_t \tilde C^\top(\tilde C \tilde V_t \tilde C^\top+ \tilde D\tilde D ^\top)^{-1} \tilde C \tilde V_t \tilde A^\top+\tilde B \tilde B^\top\nonumber
\end{align} is positive definite.
\end{lem}
\smallskip
\begin{pf}
First, notice that the block in position (2,2) of $K_{z_t|Y_{t-1}}$ is $\tilde C \tilde V_t \tilde C^\top+ \tilde D\tilde D ^\top$. The latter is positive definite because $\tilde C \tilde V_t \tilde C^\top+ \tilde D\tilde D ^\top\geq \tilde D\tilde D ^\top
 >0$. Accordingly, in order to prove that $K_{z_t|Y_{t-1}}>0$
it is sufficient to prove that $\tilde P_{t+1}$, which is the Schur complement of the block (2,2)  of $K_{z_t|Y_{t-1}}$, is positive definite. Since $\tilde V_t$ is invertible we can rewrite $\tilde P_{t+1}$ by using the Woodbury formula:
\begin{align}
 \tilde P_{t+1}=\tilde A[\tilde V_t^{-1}+\tilde C^\top(\tilde D\tilde D ^\top)^{-1} \tilde C]^{-1} \tilde A^\top+\tilde B \tilde B^\top.\nonumber
\end{align} Let $v=[\, v_1\; v_2\,]^\top\in\mathbb R^{(n+1)L}$ be such that $v_1\in\mathbb R^{n} $ and $v_2\in\mathbb R^{nL} $.  Notice that {
\begin{align}
 v^\top\tilde P_{t+1}v =\nonumber
 v^\top\tilde A[\tilde V_t^{-1}+\tilde C^\top(\tilde D\tilde D ^\top)^{-1} \tilde C]^{-1} \tilde A^\top v+v_1^\top B  B^\top v_1\nonumber
\end{align}} and thus
\begin{align}
 \label{ineq_vP_1}v^\top\tilde P_{t+1}v &\geq v_1^\top B  B^\top v_1\\
 \label{ineq_vP_2} v^\top\tilde P_{t+1}v &\geq v^\top\tilde A[\tilde V_t^{-1}+\tilde C^\top(\tilde D\tilde D ^\top)^{-1} \tilde C]^{-1} \tilde A^\top v.
\end{align}
Assume that $v^\top\tilde P_{t+1}v=0$. Since $BB^\top>0$, by (\ref{ineq_vP_1}) we have that $v_1=0$. Accordingly, the inequality in {(\ref{ineq_vP_2})} for $v_1=0$ becomes
{\small \begin{align}
 v^\top\tilde P_{t+1}v &\geq \left[\begin{array}{cc}0 &v_2^\top \end{array}\right]\tilde A[\tilde V_t^{-1}+\tilde C^\top(\tilde D\tilde D ^\top)^{-1} \tilde C]^{-1} \tilde A^\top
\left[\begin{array}{c} 0  \\ v_2  \end{array}\right].\nonumber
\end{align}}Moreover, in view of the particular structure of $\tilde A$, we have that $[\,0 \; v_2^\top\,]\tilde A=[\, v_2^\top\; 0\,]=0$ if and only if $v_2=0$. Since $[\tilde V_t^{-1}+\tilde C^\top(\tilde D\tilde D ^\top)^{-1} \tilde C]^{-1}>0$, because $\tilde V_t>0$, it follows that $v_2=0$ and thus $v=0$. We proved that if $v^\top\tilde P_{t+1} v=0$ then $v=0$, i.e. $\tilde P_{t+1}$ is positive definite.\qed
\end{pf}

Finally, if $\tilde{f}_{t}\left(\xi_{t} | Y_{t-1}\right)$ is Gaussian then, in view of (\ref{tildef}), also $\tilde{f}_{t+1}\left(\xi_{t+1} | Y_{t}\right)$ is Gaussian. Accordingly, the assumption that $\tilde{f}_{0}\left(\xi_{0}\right)$ is Gaussian, implies that $\tilde{f}_{t}\left(\xi_{t} | Y_{t-1}\right)$ is Gaussian for any $t$.

\begin{thm}\label{teo_Gauss}
Consider the state space model (\ref{ss_space}) where we recall that $\mathrm{rank} (B)=n$ and $\mathrm{rank} (D)=m$. Let $\tilde{f}_{t}(\xi_{t}|Y_{t-1} )\sim \mathcal{N}(\hat{\xi}_{t}, \tilde V_{t})$ with $\tilde V_t>0$. Then, the estimator \begin{equation}
{g}_t^0\left(y_{t}\right)= \tilde{A} \hat{\xi}_{t}+\tilde G_{t}(y_{t}- \tilde{C} \hat{\xi}_t)
\label{g0t}
\end{equation} with $$
\begin{aligned} \tilde G_{t} &=\tilde{A} \tilde V_{t} \tilde{C}^{\top }(\tilde{C} \tilde V_{t} \tilde{C}^{\top }+\tilde{D} \tilde{D}^{\top })^{-1} \end{aligned}$$
solves Problem  (\ref{minimax}). The nominal error covariance of $\xi_{t+1}$ given $y_t$ is
$$
\begin{aligned}
\tilde P_{t+1}&=\tilde{A} \tilde V_{t} \tilde{A}^{\top }-\tilde G_{t}(\tilde{C} \tilde V_{t} \tilde{C}^{\top }+\tilde{D} \tilde{D}^{\top }) {\tilde G_{t}}^{\top }+\tilde B \tilde B^{\top }
\end{aligned}
$$and the perturbed error covariance of $\xi_{t+1}$ given $Y_t$ is
$$\tilde V_{t+1}=(\tilde P_{t+1}^{-1}-\lambda^{-1} \tilde H^\top \tilde H)^{-1}.$$
If we denote $r(P)$ as the largest eigenvalue of $P$, {then} the Lagrange multiplier ${\lambda_t}>r(\tilde H \tilde P_{t+1} \tilde H^\top)$ is unique and such that
\begin{equation}
\begin{aligned} \label{gamma_t}  \gamma\left(\lambda_{t}\right)=\frac{1}{2}\left[\tr\left((I-\lambda_{t}^{-1} \tilde H^\top \tilde H \tilde P_{t+1}  )^{-1}-I\right)+\ln \operatorname{det}(I-\lambda_{t}^{-1} \tilde H^\top  \tilde H \tilde P_{t+1})\right]=c_{t}. \end{aligned}
\end{equation}
Finally, the least favorable density $\tilde f_t^0(z_t|Y_{t-1})$ corresponding to the solution of (\ref{minimax}) is a non-degenerate Gaussian density.
\end{thm}
\smallskip
\begin{pf}
As we already noticed, $\bar f_t(z_t|Y_{t-1})$ is Gaussian, and in view of Lemma \ref{lemma3}, non-degenerate. Accordingly, in view of (\ref{tildef0}), $\tilde f_t \left(z_t|Y_{t-1} \right)$ is Gaussian and non-degenerate. Let $$
\tilde{f}_{t}\left(z_{t} | Y_{t-1}\right) \sim \mathcal{N}(\tilde{m}_{z_{t} | Y_{t-1}}, \tilde{K}_{z_{t} | Y_{t-1}})$$
with $$\tilde m_{z_t|Y_{t-1}}=\left[\begin{array}{c}{\tilde{m}_{\xi_{t+1} | Y_{t-1}}} \\ {\tilde{m}_{y_t | Y_{t-1}}}\end{array}\right],$$ and $$ \tilde K_{z_t|Y_{t-1}}=\left[\begin{array}{cc}{\tilde{K}_{\xi_{t+1} | Y_{t-1}}} & {\tilde{K}_{\xi_{t+1} y_t | Y_{t-1}}} \\ {\tilde{K}_{y_t \xi_{t+1} | Y_{t-1}}} & {\tilde{K}_{y_t | Y_{t-1}}}\end{array}\right].$$
In view of (\ref{barf}) and (\ref{tildef0}), the conditional KL-divergence in (\ref{tildeE}) admits the closed-form expression
\begin{equation}
\begin{aligned} \tilde{\Es}&\left[\ln  ({\tilde{\phi}_{t}^\circ} /\phi_{t})| Y_{t-1}\right]=\tilde{\Es}\left[\ln  ( \tilde{f}_{t} /\bar{f}_{t})  | Y_{t-1}\right]\\ =&\frac{1}{2}\left[\left\|\Delta m\right\|_{K_{z_t|Y_{t-1}}^{-1}}^{2}+\tr (K_{z_t|Y_{t-1}}^{-1} \tilde{K}_{z_t|Y_{t-1}}-I )-\ln \operatorname{det} (K_{z_t|Y_{t-1}}^{-1} \tilde{K}_{z_t|Y_{t-1}} )\right]
\label{tildeE2}
\end{aligned}
\end{equation}
where $$
\Delta m=\tilde{m}_{z_t|Y_{t-1}}-m_{z_t|Y_{t-1}}
.$$ Therefore, we can rewrite the minimax game (\ref{minimax}) with respect to $\tilde f_t \left(z_t|Y_{t-1} \right)$ and $\bar f_t(z_t|Y_{t-1})$:
\begin{equation*}
(g_t^0,\tilde f^0_t )=\arg \min _{g_t \in \mathcal{\bar G}_{t}} \max_{\tilde{f}_{t} \in \mathcal{\bar B}_{t}} J_t(\tilde {f}_t,g_t)
\end{equation*}
where $$\mathcal {\bar B}_t=\left\{\tilde f_t(z_t|Y_{t-1})~s.t.~ \tilde{\Es}\left[\ln  ( \tilde{f}_{t}/\bar{f}_{t})  | Y_{t-1}\right] \leq c_t \right\},$$
$\mathcal {\bar G}_t$ is the set of estimators with finite second order moments with respect to all the densities $\tilde f_t \in \mathcal {\bar B}_t$ and
\begin{equation*}
\begin{aligned}
J_t(\tilde f_t, g_t ):=\frac{1}{2} \int_{\mathbb R^{(L+1)n+m}} \|\tilde H (\xi_{t+1}&-g_{t} (y_{t})) | Y_{t-1}\|^{2} \times \tilde f_t\left(z_t|Y_{t-1}\right)dz_t.
\end{aligned}
\end{equation*}
Next, we prove that $\tilde f^0_t$ and $g^0_t$ are such that
\begin{equation}
J_{t}(\tilde{f}_{t}, g_{t}^{0} ) \leq J_{t} (\tilde{f}_{t}^{0}, g_{t}^{0} ) \leq J_{t} (\tilde{f}_{t}^{0}, g_{t} ),
\label{Jt}
\end{equation}
where $\tilde{f}^0_{t}\left(z_{t} | Y_{t-1}\right) \sim \mathcal{N}(\tilde{m}^0_{z_{t} | Y_{t-1}}, \tilde{K}^0_{z_{t} | Y_{t-1}})$ with
\begin{equation}
\begin{aligned}
&\tilde{m}^0_{z_{t} | Y_{t-1}}=m_{z_{t} | Y_{t-1}}, \\
&\tilde{K}^0_{z_{t} | Y_{t-1}}=\left[\begin{array}{cc}{\tilde{K}_{\xi_{t+1} | Y_{t-1}}} & {K_{\xi_{t+1} y_t | Y_{t-1}}} \\ {K_{y_t \xi_{t+1} | Y_{t-1}}} & {K_{y_t | Y_{t-1}}}\end{array}\right].
\label{mK}
\end{aligned}
\end{equation}
Since $\tilde f^0_t$ is Gaussian, the optimal estimator satisfying the second inequality in (\ref{Jt}) is (\ref{g0t}). Then, it remains to prove that the least favorable density $\tilde f^0_t$ is such that (\ref{mK}) holds.

It is not difficult to see that
\begin{equation}
\begin{aligned} J_t(\tilde f_t, g_t^0 ) =&\frac{1}{2} \tr\left\{\left[\begin{array}{c}{I} \\ {-\tilde G_{t}^{\top }}\end{array}\right]\tilde H^{\top } \tilde H\left[I \quad- \tilde G_{t}\right] \times \left(\tilde{K}_{z_{t}| Y_{t-1}}+\Delta m \Delta m^{\top } \right)\right\}.
\label{Jt2}
\end{aligned}
\end{equation}
Then, based on the parametric structure of the KL divergence in (\ref{tildeE2}) and the objective function in (\ref{Jt2}), we consider the corresponding Lagrangian as a function of $\tilde m_{z_t|Y_{t-1}}$ and $\tilde K_{z_t|Y_{t-1}}$ as follows:
\begin{equation*}
\begin{aligned} \mathcal L(&\tilde{m}_{z_{t}| Y_{t-1}}, \tilde{K}_{z_{1}|Y_{t-1}}, \lambda_t)\\
=& J(\tilde{m}_{z| Y_{t-1}}, \tilde{K}_{z_t| Y_{t-1}})+\lambda_t \left(c_t-\tilde{\Es} [\ln (\tilde{f}_{t}/\bar{f}_{t})| Y_{t-1} ]\right) \\
=& \frac{1}{2} \tr\left\{\left[\begin{array}{c}{I} \\ {-{\tilde G_{t}}^{\top}}\end{array}\right]\tilde  H^{\top} \tilde H\left[I\quad-{\tilde G_{t}}\right](\tilde{K}_{z_{t}| Y_{t-1}} +\Delta m \Delta m^{\top})\right\}   \\
&-\frac{\lambda_t}{2} \Delta m^{\top} K_{z_{t}| Y_{t-1}}^{-1} \Delta m -\frac{\lambda_t}{2}   \tr(K_{z_{t}| Y_{t-1}}^{-1} \tilde{K}_{z_{t}|Y_{t-1}})+ \frac{ \lambda_t}{2} \tr(I)\\
& +\frac{\lambda_t}{2}  \operatorname{In} \operatorname{det}(K_{z_{t}| Y_{t-1}}^{-1} \tilde{K}_{z_t| Y_{t-1}})+ \lambda_t c_t\\
=& \lambda_t c_t+\frac{\lambda_t}{2}  \tr(I)+\frac{\lambda_t }{2} \ln \operatorname{det} (K_{z_{t} | Y_{t-1}}^{-1} \tilde{K}_{z_{t} | Y_{t-1}} )
\\ &+ \frac{1}{2} \tr (W(\lambda_t) \tilde{K}_{z_{t} | Y_{t-1}})+ \frac{1}{2} \Delta m^{\top} W(\lambda_t) \Delta m \end{aligned}
\end{equation*}
where $$
W(\lambda_t) \triangleq\left[\begin{array}{c}{I} \\ {-{\tilde G_{t}}^{\top}}\end{array}\right]\tilde H^\top \tilde H\left[\begin{array}{cc}{I} & {-{\tilde G_{t}}}\end{array}\right]-\lambda_t K_{z_t|Y_{t-1}}^{-1}.
$$
The first variation and the second variation of $\mathcal L$ with respect to $\tilde{m}_{z_{t}| Y_{t-1}}$ are, respectively,
\begin{equation*}
\begin{aligned}
\delta \mathcal L&(\tilde{m}_{z_{t} | Y_{t-1}}, \tilde{K}_{z_{t} | Y_{t-1}}, \lambda_t ; \delta \tilde{m}_{z_{t} | Y_{t-1}} )\\
&~~~~~=\frac{1}{2} \delta \tilde{m}^\top_{z_{t} | Y_{t-1}} W(\lambda_t) \Delta m+\Delta m^{\top} W(\lambda_t) \delta \tilde{m}_{z_{t} | Y_{t-1}},
\end{aligned}
\end{equation*}
\begin{equation*}
\begin{aligned}
&\delta^{2} {\mathcal L}(\tilde{m}_{z_{t} | Y_{t-1}},  \tilde{K}_{z_{t} | Y_{t-1}}, \lambda_t ; \delta \tilde{m}_{z_{t} | Y_{t-1}},\delta \tilde{m}_{z_{t} | Y_{t-1}})=\delta \tilde{m}_{z_{t} | Y_{t-1}} W(\lambda_t) \delta \tilde{m}_{z_{t} | Y_{t-1}}
\end{aligned}
\end{equation*}
so that $\delta^{2} {\mathcal L}<0$, for any $ \tilde{m}_{z_{t} | Y_{t-1}} \neq 0$ if and only if $W(\lambda_t)<0$, which means ${\mathcal L}$ is strictly concave if and only if $W(\lambda_t)$ is negative definite. {Next, we find the condition on $\lambda_t$ for which}  $W(\lambda_t)<0$. We denote
$M=\tilde H[\,I\; {-\tilde G_{t}}\,] O_{z_t}$, where $K_{z_{t} | Y_{t-1}}=O_{z_t} O_{z_t}^\top$, so that
\begin{equation*}
\begin{aligned} &M M^{\top} =\tilde H\left[\begin{array}{cc}{I} & { -\tilde G_{t}}\end{array}\right] O_{z_t} O_{z_t}^{\top}\left[\begin{array}{c}{I} \\ {-\tilde G^\top_{t}}\end{array}\right]\tilde H^\top \\ &=\tilde H\left(K_{\xi_{t+1}|Y_{t-1}}-K_{\xi_{t+1} y_t|Y_{t-1}} K_{y_t|Y_{t-1}}^{-1} K_{y_t \xi_{t+1}|Y_{t-1}}\right)\tilde H^\top =\tilde H \tilde P_{t+1}\tilde H^\top.\end{aligned}
\end{equation*}
Therefore, \begin{equation*}
\begin{aligned} W(\lambda_t)& =  O^ {-\top}_{z_t}\left(O_{z_t}^{\top}\left[\begin{array}{c}{I} \\ {- \tilde G_{t}^{\top}}\end{array}\right] \tilde H^{\top}
 \tilde H\left[\begin{array}{cc}{I} & -\tilde G_{t}\end{array}\right] O_{z_t}-\lambda_t I\right) O_{z_t}^{-1} \\& = O_{z_t}^{-\top}\left(M^{\top} M-\lambda I\right) O_{z_t}^{-1} \end{aligned}
\end{equation*}
and it is congruent to $\bar{W}(\lambda_t)=M^{\top} M-\lambda I$, which means $\bar W(\lambda_t)<0$ if and only if $W(\lambda_t)<0$, meanwhile, it is not difficult to see that $r(M^\top M)=r(MM^\top)=r(\tilde HP_{t+1}\tilde H^\top)$. Hence, $W(\lambda_t)<0$ as long as the Lagrange multiplier $\lambda_t>r(\tilde H \tilde P_{t+1}\tilde H^\top)$.  In such a situation the minimum is such that $\Delta m=0$, which implies
\begin{equation}\tilde{m}^0_{z_{t} | Y_{t-1}}={m}_{z_{t} | Y_{t-1}}.\label{deltam}
\end{equation}
The first variation and the second variation of $\mathcal L$ with respect to $\tilde{K}_{z_{1}|Y_{t-1}}$ are, respectively,
\begin{equation*}
\begin{aligned}\delta L& (\tilde{m}_{z_{t} |{Y_{t-1}}}, \tilde{K}_{z_{t} | Y_{t-1}}, \lambda_t ; \delta \tilde{K}_{z_{t} | Y_{t-1}} ) \\ &~~~~~=\frac{1}{2}  \tr\left\{{ \lambda_t}\tilde{K}_{z_{t} | Y_{t-1}}^{-1} \delta \tilde{K}_{z_{t} | Y_{t-1}}+W(\lambda_t) \delta \tilde{K}_{z_{t} | Y_{t-1}}\right\}\end{aligned}
\end{equation*}
and
\begin{equation*}
\begin{aligned}\delta^{2} L (\tilde{m} &_{z_{t} | Y_{t-1}}, \tilde{K}_{z_{t} | Y_{t-1}}, \lambda_t ; \delta \tilde{K}_{z_{t} | Y_{t-1}},{\delta \tilde{K}_{z_{t} | Y_{t-1}}} )  \\ &~~~~~=-\frac{1}{2} \lambda_t \tr\left\{ (\tilde{K}_{z_{t} | Y_{t-1}}^{-1} \delta \tilde{K}_{z_{t} | Y_{t-1}} )^{2}\right\}<0.\end{aligned}
\end{equation*}
Accordingly, $\mathcal L$ is strictly concave in  $\tilde{K}_{z_{t}|Y_{t-1}}$. Thus, the minimum point $\tilde{K}^0_{z_{t} | Y_{t-1}}$ is given by imposing the stationarity condition $$ \delta L(\tilde{m}_{z_{t} |{Y_{t-1}}},  \tilde{K}^0_{z_{t} | Y_{t-1}}, \lambda_t ; \delta \tilde{K}_{z_{t} | Y_{t-1}})=0.$$ The latter implies that such a point is
\begin{equation}
\tilde{K}^{0,-1}_{z_{t}| Y_{t-1}}=K^{-1}_{z_{t}| Y_{t-1}}-\frac{1}{\lambda_t}\left[\begin{array}{c}{I} \\ {-{\tilde G_{t}}^{\top}}\end{array}\right]  \tilde H^{\top} \tilde H\left[\begin{array}{cc}{I} & {-\tilde G_{t}}\end{array}\right]. \label{tildeK}
\end{equation}
Notice that the  block upper diagonal lower (UDL) form of ${K}_{z_{t}| Y_{t-1}}$ is
\begin{equation}
K_{z_{t} | Y_{t-1}}=\left[\begin{array}{cc}{I} & {\tilde G_{t}} \\ {0} & {I}\end{array}\right]\left[\begin{array}{cc}{\tilde P_{t+1}} & {0} \\ {0} & { K_{y_t | Y_{t-1}}}\end{array}\right]\left[\begin{array}{cc}{I} & {0} \\ {{ \tilde G_{t}}^\top} & {I}\end{array}\right] \label{K}
\end{equation}
and its inverse admits the following UDL decomposition
\begin{equation}
K_{z_t| Y_{t-1}}^{-1}=\left[\begin{array}{cc}{I} & {0} \\ {-{\tilde G_{t}}^{\top}} & {I}\end{array}\right]\left[ \begin{array}{cc}{\tilde P^{-1}_{t+1}} & {0} \\ {0} & { K_{y_t | Y_{t-1}}^{-1}}\end{array}\right]\left[\begin{array}{cc}{I} & {- \tilde G_{t}} \\ {0} & {I}\end{array}\right]. \label{inverseK}
\end{equation}
Therefore, substituting Equation (\ref{inverseK}) in Equation (\ref{tildeK}), we have
\begin{equation*}
\begin{aligned}
{ \tilde{K}_{z_{t} |Y_{t-1}}^{0,-1}}=\left[ \begin{array}{cc}{I} & {0} \\ {-{\tilde G_{t}}^{\top}} & {I}\end{array}\right]\left[ \begin{array}{cc}{\tilde P_{t+1}^{-1}-\lambda_t^{-1} \tilde H^{\top} \tilde H} & {0} \\ {0} & { K_{y_t | Y_{t-1}}^{-1}}\end{array}\right]\left[ \begin{array}{cc}{I} & {-{\tilde G_{t}}} \\ {0} & {I}\end{array}\right],
\end{aligned}
\end{equation*}
so that
\begin{equation}
\begin{aligned}
{ \tilde{K}_{z_{t} | Y_{t-1}}^{0}}=\left[\begin{array}{cc}{I} & {\tilde G_{t}} \\ {0} & {I}\end{array}\right]\left[\begin{array}{cc}{\tilde V_{t+1}} & {0} \\ {0} & { K_{y_t| Y_{t-1}}}\end{array}\right]\left[\begin{array}{cc}{I} & {0} \\ {{\tilde G_{t}}^{\top}} & {I}\end{array}\right]\label{tildeK2}
\end{aligned}
\end{equation}
where $$\tilde V_{t+1}=(\tilde P_{t+1}^{-1}-\lambda_t^{-1} \tilde H^{\top} \tilde H)^{-1}.$$
Then, let  $\gamma (\lambda_{t} ):=\tilde{\Es}[\ln  ( \tilde{f}_{t}^{0} / \bar{f}_{t} )  | Y_{t-1}].$
By taking into account (\ref{tildeE2}) and using Equations (\ref{deltam}), (\ref{K}) and (\ref{tildeK2}), we obtain (\ref{gamma_t}). The first derivative of $\gamma(\lambda_t)$ is
\begin{equation*}
\begin{aligned}
\frac{\partial \gamma(\lambda_t ; \delta \lambda_t)}{\partial \lambda_t}=&\frac{\lambda_{t}^{-2}}{2} \tr\left[ (I-\lambda_{t}^{-1} \tilde{H}^{\top} \tilde{H}^{\top} \tilde{P}_{t+1})^{-1}\right.\\
 &\left.\times \tilde{H}^{\top} \tilde{H}\tilde P_{t+1}(I- (I-\lambda_{t}^{-1} \tilde{H}^{\top} \tilde{H} \tilde{P}_{t+1} )^{-1})\right]\\
 =&\frac{\lambda_{t}^{-2}}{2} \tr\left[ (I-\lambda_{t}^{-1} \tilde{H}^{\top} \tilde{H}\tilde P_{t+1})^{-1} \tilde{H}^{\top} \tilde{H}\tilde P_{t+1}\right.\\
&\left. \times (I-\lambda_{t}^{-1} \tilde{H}^{\top} \tilde{H} \tilde{P}_{t+1})^{-1}((I-\lambda_{t}^{-1} \tilde{H}^{\top} \tilde{H}^{\top} \tilde{P}_{t+1})-I)\right]\\
=&-\frac{\lambda_{t}^{-3}}{2} \tr\left[((I-\lambda_{t}^{-1} \tilde{H}^{\top} \tilde{H} \tilde{P}_{t+1})^{-1} \tilde{H}^{\top} \tilde{H}\tilde P_{t+1})^{2}\right]<0.
\end{aligned}
\end{equation*}
Therefore, $\gamma(\lambda_t)$ is strictly monotone decreasing. Moreover, it is not difficult to see that \begin{equation}
\lim _{\lambda_t \rightarrow \infty} \gamma(\lambda_t)=0, \quad \lim _{\lambda_t \rightarrow r(\tilde Q)} \gamma(\lambda_t)=+\infty
\end{equation}
where $r(\tilde Q)=r(\tilde H \tilde P_{t+1} \tilde H^\top )$. As a consequence, there exists a unique Lagrangian multiplier $\lambda_t>r(\tilde Q)>0 $  such that $D(\tilde f, f)=c_t$. The fact that $\tilde f_t^0(z_t|Y_{t-1})$ is non-degenerate follows from the fact that ${\tilde K_{z_t|Y_{t-1}}^0}\geq K_{z_t|Y_{t-1}} >0$.\qed
\end{pf}

\begin{corollary} Consider the state space model (\ref{ss_space})  {where we recall:} $\mathrm{rank} (B)=n$, $\mathrm{rank} (D)=m$ and {$\xi_0\sim{\mathcal N}(\hat \xi_0,\tilde V_0)$ with  $ \tilde V_0>0$}. Then, $\tilde f_t^0(z_t|Y_{t-1})$ is non-degenerate for any $t\geq 0$.\end{corollary}
\smallskip
\begin{pf}
 We prove the claim by induction. Let $\tilde f_t(\xi_t|Y_{t-1})\sim \mathcal N(\hat \xi_t,\tilde V_t)$ with $\tilde V_t>0$. By Theorem \ref{teo_Gauss} we have that $\tilde f_t^{0}(z_t|Y_{t-1})$ is Gaussian non-degenerate and {$V_{t+1}>0$}. From, $\tilde f_t(z_t|Y_{t-1})$ we have that
$\tilde f_{t+1}(\xi_{t+1}|Y_t)\sim {\mathcal N}(\hat \xi_{t+1},\tilde V_{t+1})$ which is non-degenerate. Finally, at the initial time $t=0$, we have $\tilde f_0(\xi_0|Y_{-1}):=\tilde f_0(\xi_0)\sim \mathcal N(\hat \xi_0,\tilde V_0)$ and $\tilde  V_0$ is positive definite by assumption.\qed
\end{pf}

The resulting robust fixed-lag smoother is outlined in Algorithm \ref{code:recentEnd} where $\theta_t:=\lambda_t^{-1}$ is the risk sensitivity parameter and
\begin{align}\begin{aligned} \label{def_gamma2}
\gamma(\tilde P_{t+1},\theta_t):=\frac{1}{2}\left[\tr\left((I-\theta_t \tilde H^\top \tilde H \tilde P_{t+1}  )^{-1}-I\right)+\ln \operatorname{det}(I-\theta_t \tilde H^\top  \tilde H \tilde P_{t+1})\right].
\end{aligned}\end{align}

\begin{algorithm}[h]
{\caption{Robust { fixed-lag} smoother with {lag} L}    \label{code:recentEnd}
  \begin{algorithmic}[1]
    \Require
      $y_0 \ldots y_N$,
      $\hat \xi_{0}$,
      $\tilde V_{0}$, $\mz c_t$
      \Ensure  $\hat x_{t-L+1|t}$, $t=L-1\ldots N$
    \For{$t= 0:N $}
    \State $\tilde G_{t}=\tilde A\tilde V_{t} \tilde C^\top(\tilde C \tilde V_{t}\tilde C^\top+\tilde D\tilde D^\top)^{-1}$
          \State \label{step_algo_xi} $\hat \xi_{t+1}=\tilde A \hat \xi_{t}+\tilde G_t(y_t-\tilde C\hat \xi_t )$
      \State $\hat x_{t-L+1|t}=   \tilde H \hat \xi_{t+1}$
      \State \label{Ric_step}$\tilde P_{t+1}=\tilde A\tilde V_{t} \tilde A^\top-\tilde G_t(\tilde  C\tilde V_t \tilde C^\top+\tilde  D\tilde  D^\top){\tilde G_t^\top}+\tilde B \tilde B^\top$
 \State \label{thete_algo1}Find $\theta_t$ s.t. $\gamma(\tilde P_{t+1},\theta_t)=c_{t}$
 \State $\tilde V_{t+1}=(\tilde P^{-1}_{t+1}-\theta_t \tilde H^\top \tilde H)^{-1} $
     \EndFor
  \end{algorithmic}}
\end{algorithm}
Finally, in the case that $c_t=0$, i.e. the nominal model coincides with the actual one, it is not difficult to see that $\theta_t=0$ that is we obtain the standard fixed-lag smoother.

\begin{rmk} In the presence of a deterministic input $u_t$, then it is possible to derive the corresponding robust fixed-lag smoother by using arguments similar to the ones in \cite{RS_MPC_IET}. For instance, if the input acts only in the state equation, i.e. we have $x_{t+1}=A x_t+B v_t+u_t$, then Step \ref{step_algo_xi} in Algorithm \ref{code:recentEnd} is substituted with $\hat \xi_{t+1}=\tilde A \hat \xi_{t}+\tilde G_t(y_t-\tilde C\hat \xi_t )+w_t$, where $w_t:= [\,u_t^\top \, u_{t-1}^\top \ldots u_{t-L+1}^\top\,]^\top$.
\end{rmk}

\section{Efficient implementation}\label{sec_4}
Algorithm \ref{code:recentEnd} is not numerically robust and efficient in terms of computational burden. Since the dimension of $\tilde P_t$ and $\tilde V_t$ is proportional to $L$, which is typically large, their inversion is time consuming and not accurate. Accordingly, there is the need to develop an efficient strategy which avoids those matrix inversions as it has been done in \cite{discretetimefixedlag} for the standard fixed-lag smoother. The efficient procedure for our robust smoother is outlined in Algorithm \ref{efficent}. Next, we explain how to derive the salient steps. In what follows we always refer to the steps of Algorithm \ref{efficent} if not specified.

First, we rewrite the risk-sensitive Riccati iteration in Step \ref{Ric_step} of Algorithm \ref{code:recentEnd} as:
\begin{align}
\label{cond1}\tilde L_t&= \tilde V_t\tilde C^\top(\tilde C\tilde  V_t\tilde C^\top+\tilde  D \tilde D^\top)^{-1}\\
\label{cond2} \tilde P_{t|t} &= (I-\tilde L_t \tilde C)\tilde V_{t}\\
\label{cond3}\tilde P_{t+1}&=\tilde A \tilde P_{t|t} \tilde A^\top+\tilde B \tilde B^\top
\end{align}where $\tilde P_{t|t}:=\tilde{\Es}[(\xi_t-\hat \xi_{t|t})(\xi_t-\hat \xi_{t|t})^\top]$ and $\hat \xi_{t|t}$ is the estimator of $\xi_t$ given $Y_t$. Then, we {parameterize} $\tilde V_{t}$ and $\tilde L_t$ in blocks of $n\times n$ matrices as follows:
\begin{equation*}\begin{aligned} \tilde{V}_{t}&=\left[\begin{array}{cccccc}
V_{t} & (V_{t}^{1})^\top & \cdots & (V_{t}^{j})^\top & \cdots & (V_{t}^{L})^\top \\
V_{t}^{1} & V_{t}^{1,1} & \cdots & (V_{t}^{j, 1})^\top & \cdots &  (V_{t}^{L, 1})^\top \\
\vdots & \vdots & \ddots & \vdots & \ddots & \vdots \\
V_{t}^{j} & V_{t}^{j, 1} & \cdots & V_{t}^{j, k} & \cdots & (V_{t}^{L, k} )^\top\\
\vdots & \vdots & \ddots & \vdots & \ddots & \vdots \\
V_{t}^{L} & V_{t}^{L, 1} & \cdots & V_{t}^{L, k} & \cdots & V_{t}^{L, L}
\end{array}\right]\\
{\sy \tilde L_{t}}&{\sy =\left[\begin{array}{lllllll}(L_{t})^\top & (L^1_{t})^\top & \cdots & (L^j_{t})^\top & \cdots & (L^L_{t})^\top \end{array}\right]^\top.}\end{aligned}\end{equation*}
With some abuse of notation: $V_t$ is also denoted by $V_t^{0,0}$; $V_t^j$ is also denoted by $V_t^{j,0}$.
Substituting the above parametrizations in (\ref{cond1}), we obtain Steps \ref{sstep_4}, \ref{sstep_8}. In regard to the initial conditions: we only need $\hat x_0$ and $V_0$, we set $V_0^j=0$ and $V_0^{j,k}>0$, with $j,k>1$, are set arbitrary such that $\tilde V_0>0$, see (\ref{set_init}). Using a parametrization for $\tilde P_{t|t}$ and $\tilde P_{t+1}$ as the one for $\tilde V_t$: from (\ref{cond2}) and  (\ref{cond3}) we obtain Steps \ref{sstep_5}, \ref{sstep_7}, \ref{sstep_10} and Steps \ref{sstep_14}, \ref{sstep_16}, \ref{sstep_18}, respectively. {\sy Regarding Step \ref{step_12}, recall that  from Step \ref{step_algo_xi} in Algorithm \ref{code:recentEnd} we have
\begin{align}
\label{asterisco}\hat \xi_{t+1}=\tilde A \hat \xi_{t}+\tilde G_t (y_t-\tilde C\hat \xi_t ),
\end{align}
where $\tilde G_t = \tilde A \tilde L_t.$ Notice that $\hat \xi_{t}$ is the predictor of $\xi_{t}$ given $Y_{t-1}$, and it can be partitioned as $\hat \xi_{t}:=[\,( \hat x_{t|t-1}^0)^{\top}\; (\hat x^1_{t|t-1})^{\top} \; \cdots \; (\hat x^j_{t|t-1})^{\top} \; \cdots \; (\hat x^L_{t|t-1})^{\top}\,]^{\top}$
where $\hat x^j_{t|t-1} = \hat x_{t-j|t-1}, j \geq 0$. Substituting the definitions of  $\tilde A$, $\tilde C$, $\tilde L_t$ in (\ref{asterisco}) we obtain Step \ref{step_12}.} It remains to
find an efficient way to find the risk sensitivity parameter $\theta_t$. Indeed, in order to evaluate  $\gamma(\tilde P_{t+1},\cdot)$ for a specific $\theta$ in Step \ref{thete_algo1} of Algorithm \ref{code:recentEnd}, we have to perform the inversion and the eigenvalue decomposition of a  matrix whose dimension is proportional to $L$. The next result shows that it is possible to find $\theta_t$ by considering a function which requires to perform the inversion and the eigenvalue decomposition of matrices of dimension $n\times n$, see Step \ref{sstep_21}. In this way, $\theta_t$ can be computed in a numerically robust way in the case that $L$ is large.

\begin{proposition}
Consider  $\gamma$ defined in (\ref{def_gamma2}). Then, $\theta_t$ is the unique solution to $\gamma(\tilde P_{t+1},\theta_t)=c_t$ if and only if $\theta_t$ is the unique solution to  $\bar \gamma(P^{L,L}_{t+1},\theta_t)=c_t$ where
$$\begin{aligned}
\bar \gamma(P^{L,L}_{t+1}, \theta_t )=&-\frac{1}{2}\left\{\tr\left[ P_{t+1}^{L,L} (P^{L,L}_{t+1}-\theta_t I )^{-1} \right]\right.\\
&\left.+\ln \operatorname{det}\left[I-P^{L,L}_{t+1} {(P^{L,L}_{t+1}-\theta_t I)^{-1}}\right]\right\}=c_{t}. \end{aligned}$$
\end{proposition}

\begin{pf}
First, notice that $\tilde H \tilde P_{t+1} \tilde H^\top=P_{t+1}^{L,L}$. Thus, condition
$$\theta_t < r (\tilde H \tilde P_{t+1} \tilde H^\top )^{-1}$$ is equivalent to $\theta_t<r(P_{t+1}^{L,L})^{-1}$.  By (\ref{def_gamma2}), we know
\begin{equation} \label{ris_gamma0} \gamma (\tilde P_{t+1},\theta_{t})=\frac{1}{2}\left[\tr(K_1)+\ln \operatorname{det}(K_2)\right]
\end{equation}
where $$K_1=(I-\theta_{t} \tilde H^\top \tilde H \tilde P_{t+1} )^{-1}-I, \; \; K_2=I-\theta_t \tilde H^\top  \tilde H \tilde P_{t+1}.$$ Then, \begin{align}\label{ris_gamma1}
&\tr (K_1)=\tr \left\{\left[(\tilde P^{-1}_{t+1}-\theta_t \tilde H^\top \tilde H )\tilde P_{t+1}\right]^{-1}\right\}-(L+1)n\nn\\
&=\tr\left[\tilde P^{-1}_{t+1}(\tilde P^{-1}_{t+1}-\theta_t \tilde H^\top \tilde H)^{-1}\right]-(L+1)n\nn\\
&=\tr\left\{\tilde P^{-1}_{t+1}\left[\tilde P_{t+1}-\tilde P_{t+1}\tilde H^\top \times (\tilde H \tilde P_{t+1}  \tilde H^\top -\theta_t^{-1}I)^{-1} \tilde H \tilde P_{t+1} \right] \right\}-(L+1)n \nn\\
&=\tr\left\{\left[I- (\tilde H \tilde P_{t+1}  \tilde H^\top -\theta_t^{-1}I)^{-1}   \tilde H \tilde P_{t+1} \tilde H^\top \right] \right\}-n\nn \\
&=-\tr (\Gamma_{t+1}^{-1}P^{L,L}_{t+1})
\end{align}
where $\Gamma_{t+1}:=P_{t+1 | t}^{L, L}-\theta_{t}^{-1} I_{n}$ and we exploited the Woodbury matrix identity. Then,
\begin{align}\label{ris_gamma2}
\ln \operatorname{det} (K_2)&=\ln \operatorname{det}\left[(\tilde P^{-1}_{t+1}-\theta_t \tilde H^\top \tilde H ) \tilde P_{t+1}\right]\nn\\
&=-\ln \operatorname{det}\left[\tilde P^{-1}_{t+1}(\tilde P^{-1}_{t+1}-\theta_t \tilde H^\top \tilde H )^{-1}\right]\nn\\
&=-\ln \operatorname{det}\left[\tilde P^{-1}_{t+1}(
\tilde P_{t+1} - \tilde P_{t+1} \tilde H ^\top   \Gamma_{t+1}^{-1} \tilde H  \tilde P_{t+1})\right]\nn\\
&=-\ln \operatorname{det}(I-\tilde H^\top\Gamma_{t+1}^{-1}\tilde H \tilde P_{t+1} )\nn\\
&=-\ln \operatorname{det}(I-\tilde H \tilde P_{t+1} \tilde H^\top \Gamma_{t+1}^{-1})\nn\\
&=-\ln \operatorname{det}(I-P^{L,L}_{t+1}\Gamma_{t+1}^{-1})
\end{align}
where we exploited the fact that $\operatorname{det}(I_m+AB)=\operatorname{det}(I_n+BA)$ and $A \in \mathbb R^{m\times n}$, $B \in \mathbb R^{n\times m}$. Substituting (\ref{ris_gamma1})
and (\ref{ris_gamma2}) in (\ref{ris_gamma0}), we get the claim.\qed
\end{pf}
Using the Woodbury formula, we have
\begin{align*}
\tilde V_{t+1}&=(\tilde P_{t+1}^{-1}-\theta_t \tilde H^\top\tilde H )^{-1}\\
&= \tilde P_{t+1}- \tilde P_{t+1}\tilde H^\top  \Gamma_{t+1}^{-1}\tilde H^\top \tilde P_{t+1}
\end{align*}
and using the parametrization of $\tilde V_{t+1}$ and $\tilde P_{t+1}$ in terms $V_{t+1}^{j,k}$ and $P_{t+1}^{j,k}$, we obtain Steps \ref{sstep_22}, \ref{sstep_23}, \ref{sstep_25}, \ref{sstep_27}.
{\sy Finally, Step 31 is derived as follows. The estimate of $\xi_{t+1}$ given $Y_{t+1}$ is obtained by the update step of the standard Kalman filter: \begin{align}
\label{triangolo}\hat \xi_{t+1|t+1} = \tilde A \hat \xi_{t|t} + \tilde L_{t+1} (y_{t+1} -\tilde C \tilde A \hat \xi_{t|t}).\end{align}
Notice that $\hat \xi_{t|t}$ can be partitioned as $$ \hat \xi_{t|t}:=[\,(\hat x_{t|t}^0)^{\top} \; (\hat x^1_{t|t})^{\top} \; \cdots \; (\hat x^j_{t|t})^{\top} \; \cdots \; (\hat x^L_{t|t})^{\top}\,]^{\top}$$where  $ \hat x^{j}_{t|t} = \hat x_{t-j|t}, j \geq 0$.
Substituting the definitions of $\tilde A$, $\tilde C$, $\tilde L_t$ in (\ref{triangolo}), we obtain
\begin{align}\label{cerchio2} \hat x^j_{t+1|t+1}=\hat x^{j-1}_{t|t} +  L^j_{t+1} (y_{t+1}-C A \hat x_{t|t}), \quad{ j\geq 0}.\end{align}
Note that when $ j=0$, $\hat x^{j-1}_{t|t} = \hat x^{-1}_{t|t} = \hat x_{t+1|t} = A \hat x_{t|t}.$ In addition, from Equation (29), it is not difficult to see  (the derivation is the same of the one in \cite{discretetimefixedlag})  $$ L^j_{t+1} = V^j_{t+1} C^{\top} (C V_{t} C^\top + DD^\top)^{-1} ={ V^j_{t+1} V^{-1} _{t+1} L_{t+1}}, \quad j \geq 0$$where $L^0_{t+1}=L_{t+1}$ and while $V^0_{t+1}=V_{t+1}$.
Substituting the latter in (\ref{cerchio2}) we have \begin{align}\label{trapezio} \hat x^j_{t+1|t+1} = \hat x^{j-1}_{t|t} +  V^j_{t+1} V_{t+1}^{-1} L_{t+1} (y_{t+1}-C A \hat x_{t|t}).\end{align}
from which   Step 31  can be established.}

 \begin{algorithm}[H]\scriptsize
  \caption{Robust efficient version }    \label{efficent}
  \begin{algorithmic}[1]
    \Require
      $y_0 \ldots y_N$,
      $\hat x_{0}$,
    $V_{0}$,
      $\mz c_t$.
    \Ensure  $\hat x_{t-L+1|t}, t=L-1\ldots N$
      \State Set $V_{0}^{j,k}=0$ for $j\neq k$, $V_{0}^{j,j}$ positive definite  for $j>0$
        \For{$t=0:N$}
                   \State  \label{sstep_4}$L_{t}=V_{t}C^\top(CV_{t}C^\top+DD^\top)^{-1}$
                        \State \label{sstep_5} $P_{t|t}=V_{t}(I-C^\top L_{t}^\top)$;
                \For{$j=1:L$}
                    \State \label{sstep_7} $P^j_{t|t}=V^j_{t}(I-C^\top L_{t}^\top)$;
                    \State \label{sstep_8} $L^j_{t}=V^j_{t}C^\top(CV_{t}C^\top+R)^{-1}$;
                        \For{$k=1:L$}
                            \State \label{sstep_10} $P^{j,k}_{t|t}=V^{j,k}_{t}-V^j_{t}C^\top(L^k_{t})^\top$;
                        \EndFor
                \EndFor
        \State \label{step_12}$\hat x_{t+1}=A \hat x_{t}+AL_{t}(y_{t}-C \hat x_{t})$;
       \State \label{sstep_14} $P_{t+1}=AP_{t|t}A^\top+BB^\top$;
                \For{$j=1:L$}
                    \State \label{sstep_16} $P^j_{t+1}=P^{j-1}_{t|t}A^\top$;
                        \For{$k=1:L$}
                            \State \label{sstep_18} $P^{j,k}_{t+1}=P^{j-1,k-1}_{t|t}$;
                        \EndFor
                \EndFor
            \State \label{sstep_21}Find $\theta_t $ s.t. ${ \bar \gamma}(P^{L,L}_{t+1},\theta_t)=c_t$;
            \State \label{sstep_22} $\Gamma_{t+1}=P^{L,L}_{t+1}-\theta^{-1}_t I_n$;
            \State \label{sstep_23} $V_{t+1}=P_{t+1}-(P^L_{t+1})^\top \Gamma_{t+1}^{-1} P^L_{t+1}$;
                \For{$j=1:L$}
                    \State \label{sstep_25} $V^j_{t+1}=P^j_{t+1}-(P^{L,j}_{t+1})^\top\Gamma_{t+1}^{-1} P^L_{t+1}$;
                        \For{$k=1:L$}
                            \State \label{sstep_27} $V^{j,k}_{t+1}=P^{j,k}_{t+1}-(P^{L,j}_{t+1})^\top\Gamma_{t+1}^{-1} P^{L,k}_{t+1}$;
                        \EndFor
                \EndFor
          \EndFor
   {  \For{$t= L-1:N$} \vspace{0.5cm}
    \State \label{sstep_33} {\small \vspace{-0.85cm} \begin{align}\hat x_{t-L+1|t}=\hat x_{t-L+1}
      +\hspace{-0.4cm}\sum_{j=t-L+1}^{t} V_{j}^{j+L-t-1} V_{j}^{-1}L_j \left(y_j-C \hat{x}_{j}\right)\nn\end{align}}
    \EndFor }
  \end{algorithmic}
\end{algorithm}

 \begin{figure}[t]
    \centering
    \includegraphics[width=0.65\textwidth]{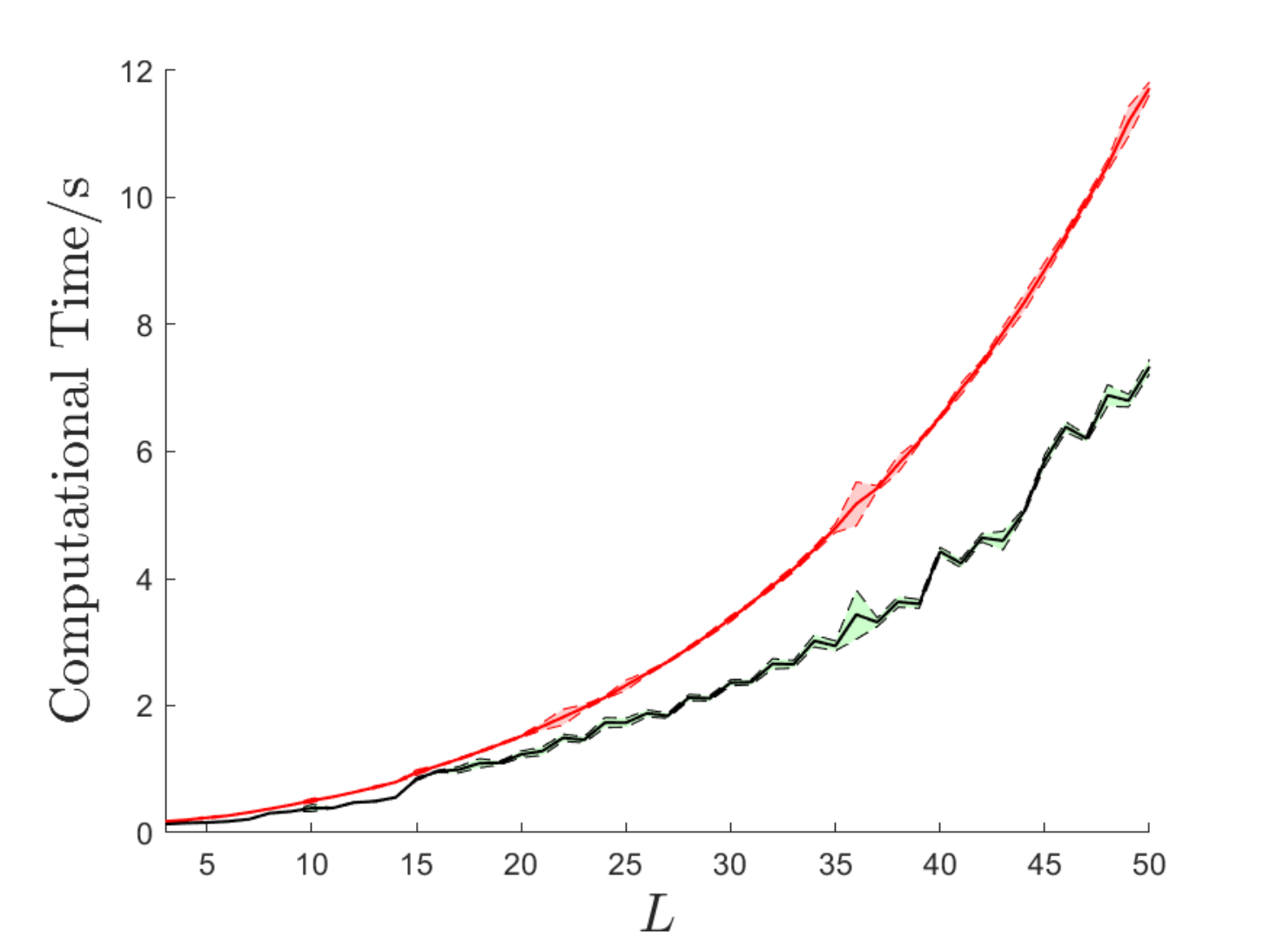}
    \caption{\footnotesize{Computational time of ARFLS (corresponding to the red part) and RFLS (corresponding to the green part).}}
    \label{time}
\end{figure}

\textbf{Computational complexity.} {\sy  We perform
the asymptotic analysis of the computational complexity, understood as the number of floating point operations (flops) by using big O notation, of the efficient version of the robust fixed-lag smoother (RFLS), i.e. Algorithm \ref{efficent}, versus the one of the augmented robust fixed-lag smoother (ARFLS), i.e. Algorithm \ref{code:recentEnd}.

First, referring to ARFLS, Steps 2-5 have the same complexity of the standard Kalman predictor, that is:
$$O(n^3(L+1)^3) + O(mn^2(L+1)^2) + Q(m^2n(L+1)) + O(m^3). $$
Here, it is worth noting that $\tilde Q=\tilde B \tilde B^{\top}$ and $\tilde R=\tilde D \tilde D^{\top}$ are computed offline. {\sy Then, in regard to Step 6, the complexity to evaluate $\gamma(\tilde P_{t+1}, \theta_t)$ is $O(n^3(L+1)^3)$, see  \cite[Section 13.1 and 13.4]{skiena1998algorithm}. Then, the computation of $\theta_t \in (0, r(\tilde P_{t+1})^{-1})$ is accomplished by a bisection method, see Algorithm 2 in \cite{zenere2018coupling}.  Since at each step we spend constant time to reduce the problem to an instance half its size  \cite[Section 4.10.2]{skiena1998algorithm}}, the complexity of Step 6 is $$O( n^3(L+1)^3 \textrm{log}_2(r(\tilde P_{t+1})^{-1}/\varepsilon) )$$
where $\varepsilon>0$ is the selected accuracy, i.e. the solution found satisfies the condition $|\gamma(\tilde P_{t+1}, \theta_t) - c_t| \leq \varepsilon$.  Step 7 has complexity $O( n^3(L+1)^3 )$.  Thus, the computational complexity of Algorithm 1 is:
\begin{align*}  O(n^3&(L+1)^3) + O(mn^2(L+1)^2) + Q(m^2n(L+1)) + O(m^3) \\ &+  O( n^3(L+1)^3   \textrm{log}_2(r(\tilde P_{t+1})^{-1}/\varepsilon) ).\end{align*} Accordingly, the complexity of Algorithm 1 with respect to the instance $L$ is $O(L^3)$.

Referring to RFLS, the total complexity of Steps 3-4, 12-13 and 21-22 is $O(n^3) + O(mn^2) + Q(m^2n) + O(m^3)$. Then, the complexity of Steps 5-11 is $O(n^3L) + O(mn^2L^2) + Q(m^2nL^2) + O(Lm^3)$; the complexity of Steps 14-19 is $O(n^3L)$; the complexity of Steps 23-28 is $O(n^3 L^2)$. Next, Step 20 has the complexity of $O( n^3 \textrm{log}_2(r(P^{L,L}_{t+1})^{-1}/\varepsilon) )$; it is worth noting thay in Step 20 the computation of $\theta_t$ is done by using the same bisection method of Step 6 in Algorithm 1: the difference is the dimension of matrices $\tilde P_{t+1}$, of dimension $n(L+1)$, and $P^{L,L}_{t+1}$, of dimension $n$. Finally, Step 31 has the complexity of $O(mn^2L)+O(n^3L).$
Hence, the computational complexity of Algorithm 2 is:
$$ O(n^3L^2) + O(mn^2L^2) + Q(m^2nL^2) + O(m^3L) +  O( n^3  \textrm{log}_2(r(P^{L,L}_{t+1})^{-1}/\varepsilon) ). $$
Thus, the complexity of Algorithm 2 with respect to the instance  $L$ is $O(L^2)$. We conclude that Algorithm 2 is computationally more efficient    than Algorithm 1 and this advantage will become more pronounced as $L$ grows.

Finally, we also analyze the computational time with respect to the lag $L$ through a Monte Carlo study. In the latter, the lag ranges from $L = 3$ up to $L = 50$. Each case is composed by 100 trials. In each trial, the matrices $A,B,C,D$ of Model (1) with $n=2$ and $m=1$ are randomly generated as follows. Each entry is drawn according to a uniform distribution in the interval $[0, 1]$.
Then, matrix $A$ is rescaled in such a way that its maximum eigenvalue (in modulus) is equal to 0.95.
Then, an output sequence $Y_N$ with $N = 500$ is generated.
Fig. 1 shows the average value of the computational time over 100 trials required by RFLS
and ARFLS with $c_t=10^{-3}$  to estimate the state trajectory from
$Y_N$.
The results were obtained using a Huawei MateBook X Pro Laptop
with Intel Core I5-8250U CPU and 8GB RAM. The dashed lines defines the
corresponding confidence intervals (with level 0.95).
It is possible to note
that the computational time of these two algorithms grows  polynomially.
As expected, the growth rate of RFLS is much smaller than  the one of ARFLS,  i.e. RFLS drastically reduces the computational time.}

\section{Least-Favorable Model} \label{sec_5}
In order to evaluate the performance of the robust fixed-lag smoother, we need to construct its least favorable model solution to (\ref{minimax}). 
The latter can be characterized over a finite time interval $[0,T]$ by using arguments similar to the ones in \cite[Section V]{yi2020lowTAC}.  More precisely, the least favorable model takes the form
\begin{equation}\begin{aligned} \label{least_favorable_model}
\eta_{t+1} &=\bar{A}_{t} \eta_{t}+\bar{B}_{t} \epsilon_{t} \\
y_{t} &=\bar{C}_{t} \eta_{t}+\bar{D}_{t} \epsilon_{t}
\end{aligned}\end{equation}
where $\eta_{t} \triangleq [\,\xi_{t}^\top\;
\tilde e_{t}^\top\,]^\top$. Moreover,
\begin{align}\label{matrices_LF}\begin{aligned}
&\bar{A}_{t}:=\left[\begin{array}{cc}
\tilde A & \tilde B \tilde S_{t} \\
0 & \tilde A-\tilde G_{t} \tilde C+(\tilde B-\tilde G_{t} \tilde D) \tilde S_{t}
\end{array}\right]\\
&\bar{B}_{t}:=\left[\begin{array}{c}
\tilde B \\
\tilde B-\tilde G_{t} \tilde D
\end{array}\right]L_t\\
&\bar{C}_{t}:=\left[\begin{array}{lll}
\tilde C &\tilde D \tilde S_{t}
\end{array}\right], \quad \bar{D}_{t}:=\tilde D L_{t}.\\
\end{aligned}\end{align}  The matrices above are computed through the backward recursion illustrated in Algorithm \ref{code:least}.

\begin{algorithm}[H]
{\caption{ Backward recursion}   \label{code:least}
  \begin{algorithmic}[1]
    \Require
      $\tilde G_0 \ldots \tilde G_N$, $\theta_0\ldots \theta_N$,
      $\tilde \Omega_{N+1}^{-1}$
      \Ensure  $\bar A_t, \bar B_t, \bar C_t, \bar D_t, t=0 \ldots N $
\State     $\tilde \Omega^{-1}_{N+1}=0$
    \For{$t= N:0 $}
      \State $W^{-1}_{t+1}=\tilde \Omega_{t+1}^{-1}+\theta_t\tilde H^\top \tilde H $
      \State $\tilde{K}_{v_{t}}=(I-(\tilde B-\tilde G_{t} \tilde D)^{\top} W_{t+1}^{-1}(\tilde B-\tilde G_{t} \tilde D))^{-1}$
      \State $\tilde S_{t}=\tilde{K}_{v_{t}}(\tilde B-\tilde G_{t} \tilde D)^{\top } W_{t+1}^{-1}(\tilde A-\tilde G_{t} \tilde C)$
      \State Compute $L_t$ such that $\tilde{K}_{v_{t}}=L_t L_t^\top$
      \State Compute $\bar{A}_{t}$, $\bar{B}_{t}$, $\bar{C}_{t}$, $\bar{D}_{t}$ as in (\ref{matrices_LF})
     \State $\tilde \Omega_{t}^{-1}=(\tilde A-\tilde G_{t} \tilde C)^{\top } W_{t+1}^{-1}(\tilde A-\tilde G_{t} \tilde C)+\tilde S_{t}^{\top } \tilde{K}_{v_{t}}^{-1} \tilde S_{t}$
     \EndFor
  \end{algorithmic}}
\end{algorithm}

It is worth noting that (\ref{least_favorable_model}) is the least favorable model corresponding to the augmented state $\xi_t$. It is then natural to wonder whether such a least favorable model reduces to a least favorable model corresponding to the state $x_t$. The answer is affirmative. This justifies why in the minimax problem (\ref{minimax}) we did not need to impose that $\tilde \phi_t(z_t|x_t)$ preserve the same structure of the augmented state space model in (\ref{model}). Substituting $\bar A_t$, $\bar B_t$, $\bar C_t$, $\bar D_t$ in (\ref{least_favorable_model}), we obtain
\begin{equation*}\begin{aligned}
\xi_{t+1} &=\tilde{A} \xi_{t}+\tilde{B} (\tilde S_t\tilde e_t+L_t\epsilon_{t}) \\
y_{t} &=\tilde{C} \xi_{t}+\tilde{D} (\tilde S_t\tilde e_t+L_t\epsilon_{t})
\end{aligned}\end{equation*} which is the augmented state space of the least favorable model
\begin{equation}\label{LF_x}\begin{aligned}
x_{t+1} &={A} x_{t}+{B} (\tilde S_t\tilde e_t+L_t\epsilon_{t}) \\
y_{t} &={C} x_{t}+{{D}} (\tilde S_t\tilde e_t+L_t\epsilon_{t}).
\end{aligned}\end{equation}
Finally, consider a fixed-lag smoother of the form
\begin{align*}
\hat x_{t-L+1|t}^\prime=\hat x_{t-L+1}^\prime+\sum_{j=t-L+1}^t G_j^{\prime \,L-t+j} (y_j-C\hat x_{j}^\prime)
\end{align*} where $G_{t}^{\prime\, j}$, with $0\leq j\leq L$, are the arbitrary gains of a fixed-lag smoothing algorithm.
It is not difficult to see that we can rewrite such a smoother as $\hat x_{t-L+1|t}^\prime = \tilde H\hat \xi_{t+1}^\prime$ where
\begin{align*}\hat{\xi}_{t+1}^{\prime}&=\tilde A \hat{\xi}_{t}^{\prime}+\tilde G_{t}^{\prime}(y_{t}-\tilde C \hat \xi_t^{\prime} )\\
\tilde G_t^{\prime}&= \left[\begin{array}{ccc} (G_t^{\prime\, 0})^\top & \ldots  &  (G_t^{\prime\, L})^\top  \end{array}\right]^\top.
\end{align*}
To evaluate its performance under the least favorable model in (\ref{LF_x}), we define the corresponding smoothing error $e_{t-L+1}^{\prime}= x_{t-L+1}-\hat x_{t-L+1|t}^{\prime}$. Then it is not difficult to see that $e_{t-L+1}^\prime$ is zero mean and with covariance matrix
\begin{equation*}
\bar \Pi_{t+1}=\tilde H \tilde \Pi_{t+1} \tilde H^\top
\end{equation*}
where: $\tilde \Pi_{t+1}$ is the $(L+1)n\times (L+1)n$ submatrix of $\Pi_{t+1}$ in position $(1,1)$; $\Pi_{t+1}$ is the solution to the Lyapunov equation \begin{equation*}\begin{aligned}
\Pi_{t+1}&=\left(\bar{A}_{t}-\left[\begin{array}{c}
\tilde G_{t}^{\prime} \\
0
\end{array}\right] \bar{C}_{t}\right) \Pi_{t}\left(\bar{A}_{t}-\left[\begin{array}{c}
\tilde G_{t}^{\prime} \\
0
\end{array}\right] \bar{C}_{t}\right)^{\top } \\
&\hspace{0.5cm}+\left(\bar{B}_{t}-\left[\begin{array}{c}
\tilde G_{t}^{\prime} \\
0
\end{array}\right] \bar{D}_{t}\right)\left(\bar{B}_{t}-\left[\begin{array}{c}
\tilde G_{t}^{\prime} \\
0
\end{array}\right] \bar{D}_{t}\right)^{\top }
\end{aligned}\end{equation*}
with initial condition  $\Pi_0=\mathbbm{1}_2\otimes
\tilde V_0$ and $\mathbbm{1}_2$ is the $2\times 2$ matrix whose entries are equal to one.

\section{Robust Target Tracking}\label{sec_7}
We compare the performance of the robust and standard fixed-lag smoothers in a maneuvering target tracking problem where model uncertainty is present. More precisely, we consider as  nominal model the second-order Singer model with an exponentially autocorrelated noise, see \cite{farahmand2011doubly,4103555,yi2017online}. The state vector is defined as $x:=\left[p^{lat},  v^{lat}, p^{lon}, v^{lon}\right]^{\top}$ where $p^{lat}$ and $v^{lat}$ denote the target position and velocity along the latitudinal direction, respectively;  $p^{lon}$ and $v^{lon}$ denote the target position and velocity along the longitudinal direction, respectively. This model can be written as  (\ref{ss_space}) with
$$
A=\left[\begin{array}{cccc}
1 & T & 0 & 0 \\
0 & 1 & 0 & 0 \\
0 & 0 & 1 & T \\
0 & 0 & 0 & 1
\end{array}\right]
$$
where $T=0.01$ denotes the sampling period;  $B$ is such that
$BB^{\top}=Q$ and
$$
Q=2 \alpha \sigma_{m}^{2}\left[\begin{array}{cccc}
T^3/3  & T^2/2 & 0 & 0 \\
T^2/2  & T     & 0 & 0 \\
0      & 0     & T^3/3 & T^2/2 \\
0      & 0     & T^2/2  & T
\end{array}\right].
$$
where we assume that $\sigma_{m}^{2}=5$ and $1/\alpha=0.5$ which are the instantaneous variance of the velocity and the time constant of the target velocity autocorrelation, respectively. Moreover, $$C=\left[\begin{array}{cccc}1 & 0 & 0 & 0 \\ 0 & 0 & 1 & 0\end{array}\right],\quad DD^{\top}= I_2$$  and the output $y_t$ denotes the noisy position measurements along the two directions. Finally, $x_0$ is Gaussian distributed with zero mean and covariance matrix $V_0=\operatorname{diag}(50,5,50,5)$.

In practice, the nominal model above does not coincide with the actual one (e.g. the nominal parameters $\alpha$ and $\sigma_{m}^{2}$ are typically imprecise). In what follows, we assume that the actual model belongs to the ambiguity set $\mathcal B_t$. More precisely, we consider two scenarios: the first one considers the ambiguity set with $c_t=10^{-3}$, while the second one with $c_t=5\cdot 10^{-3}$, i.e. the latter is larger than the former. We compare the proposed robust fixed-lag smoother, denoted by RFLS, and the standard fixed-lag smoother, denoted by FLS, both with $L=20$.

\begin{figure}[htb]
\vspace{-0.1in}
\centering
\subfigure[Ambiguity set with $c=10^{-3}$]{\includegraphics[width=2.5in]{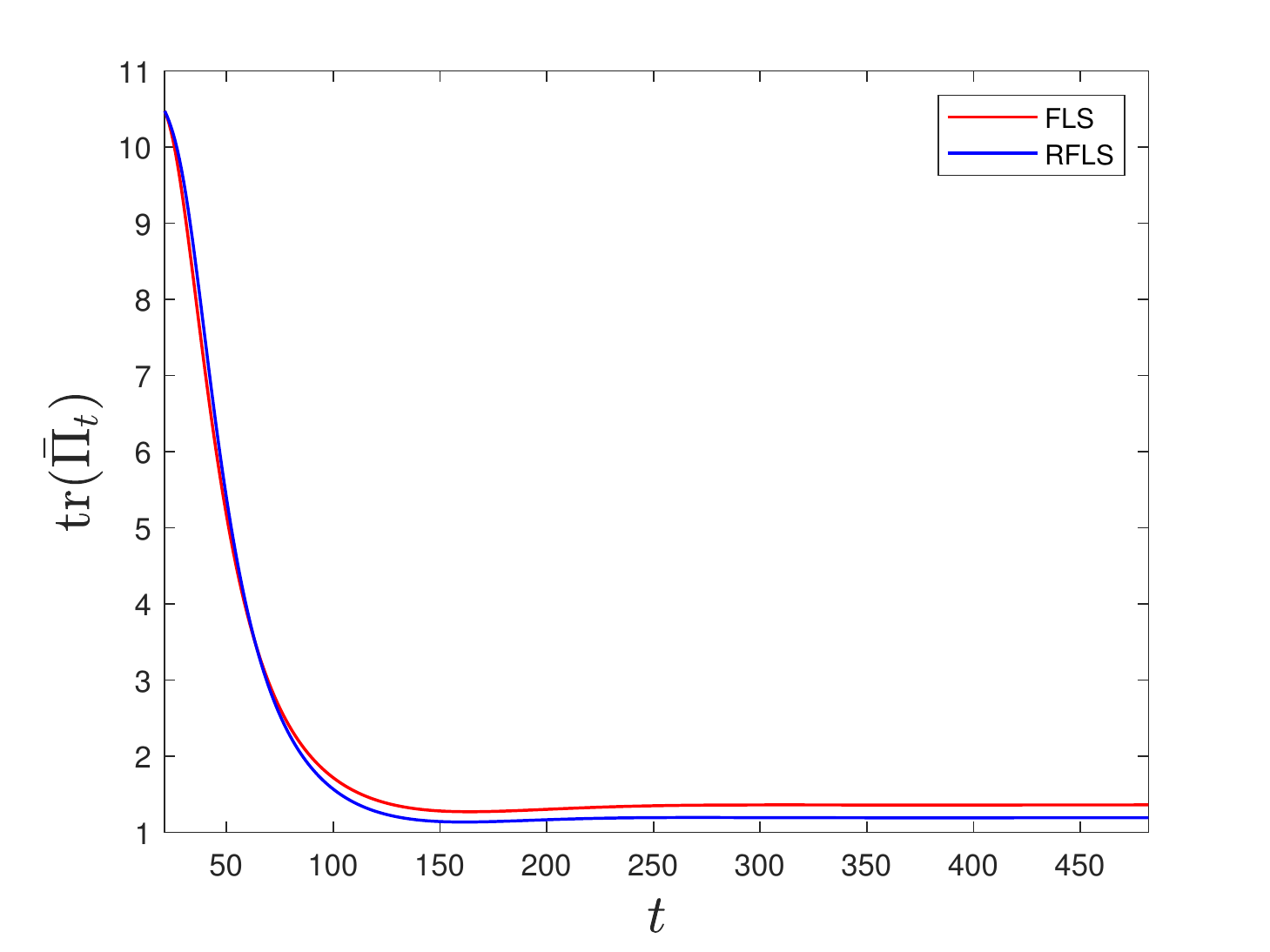}\label{track_1}}\ \
\hspace{3pt}
\subfigure[Ambiguity set with $c=5\cdot10^{-3}$]{\includegraphics[width=2.5in]{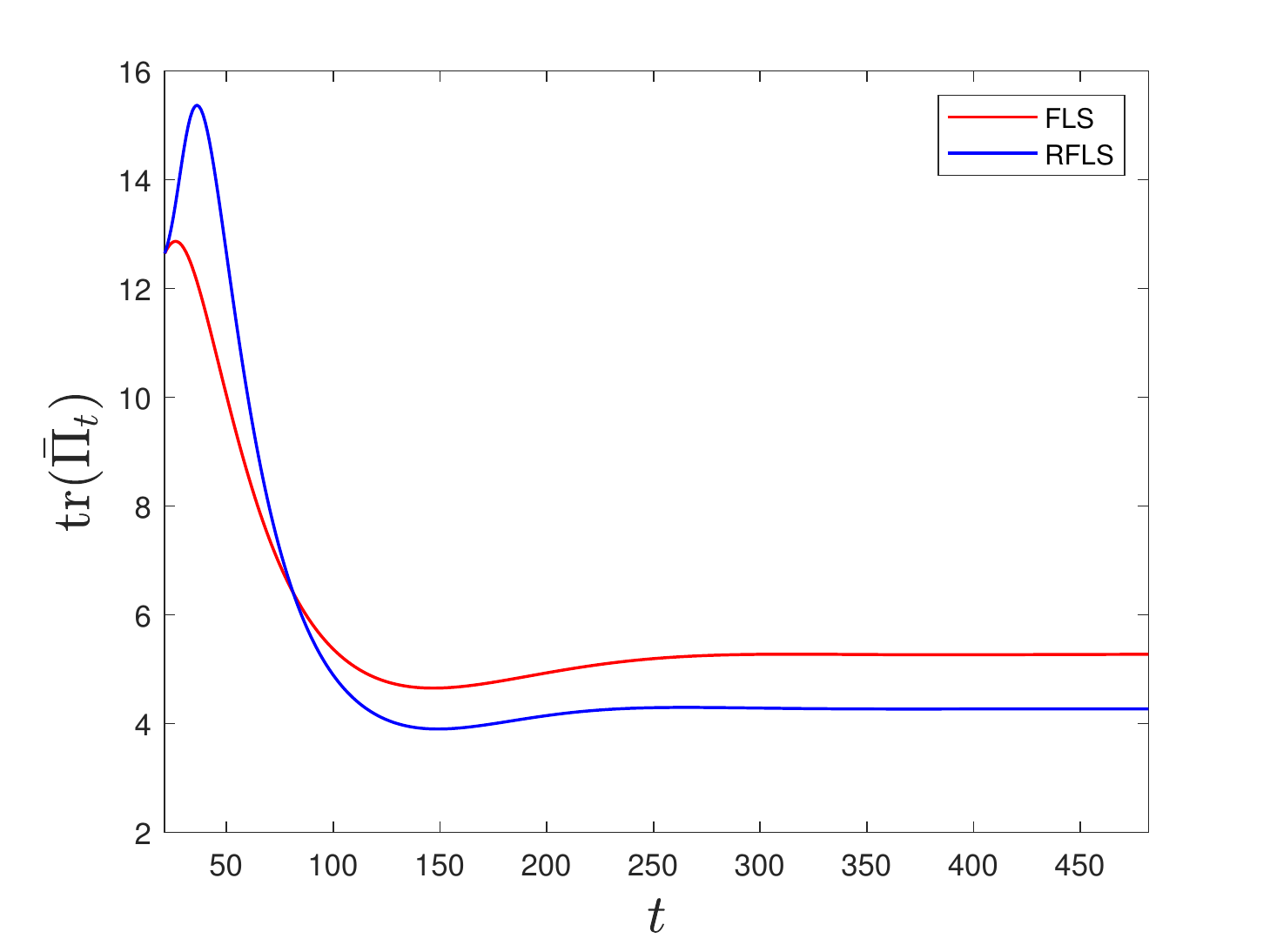}\label{track_2}}\ \
\vspace{-0.05in}
\caption{ \footnotesize{Scalar variances of the smoothing error for RFLS and FLS under the least favorable models in the different ambiguity sets.}}     \label{track}
\end{figure}

Fig. \ref{track_1} shows the variances of the smoothing error under the least favorable model in the first ambiguity set with $c_t=10^{-3}$, while Fig.  \ref{track_2} shows the variances of the smoothing error under the least favorable model in the second ambiguity set with $c_t=5\cdot10^{-3}$. It is possible to see that RFLS outperforms FLS. It is also worth noting that the higher $c_t$ is, the more RFLS outperforms FLS.
%

\begin{figure}[htb]
\centering
\subfigure[Ambiguity set with $c=10^{-3}$]{\includegraphics[width=2.5in]{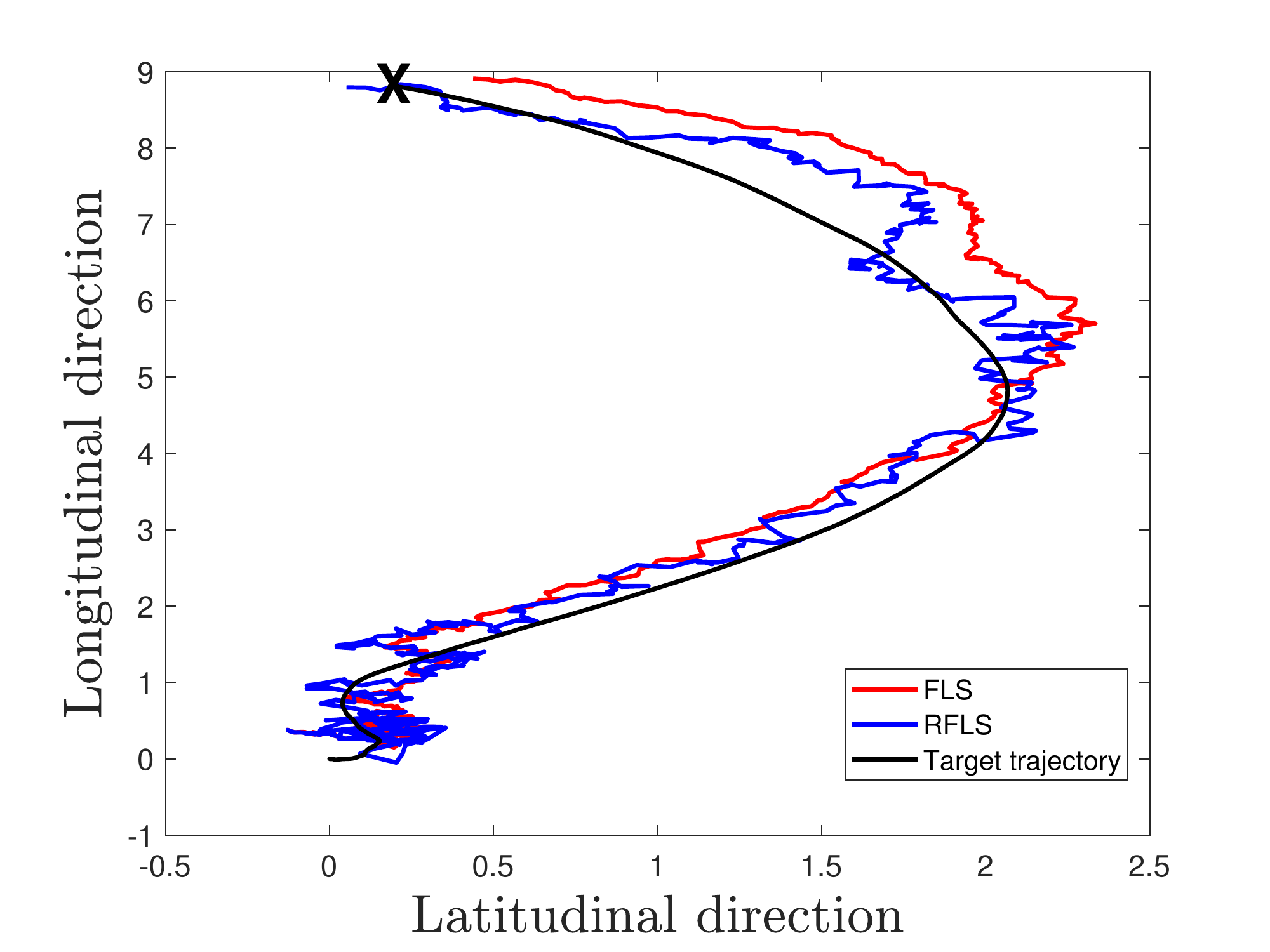}\label{pos_1}}\ \
\hspace{3pt}
\subfigure[Ambiguity set with $c=5\cdot10^{-3}$]{\includegraphics[width=2.5in]{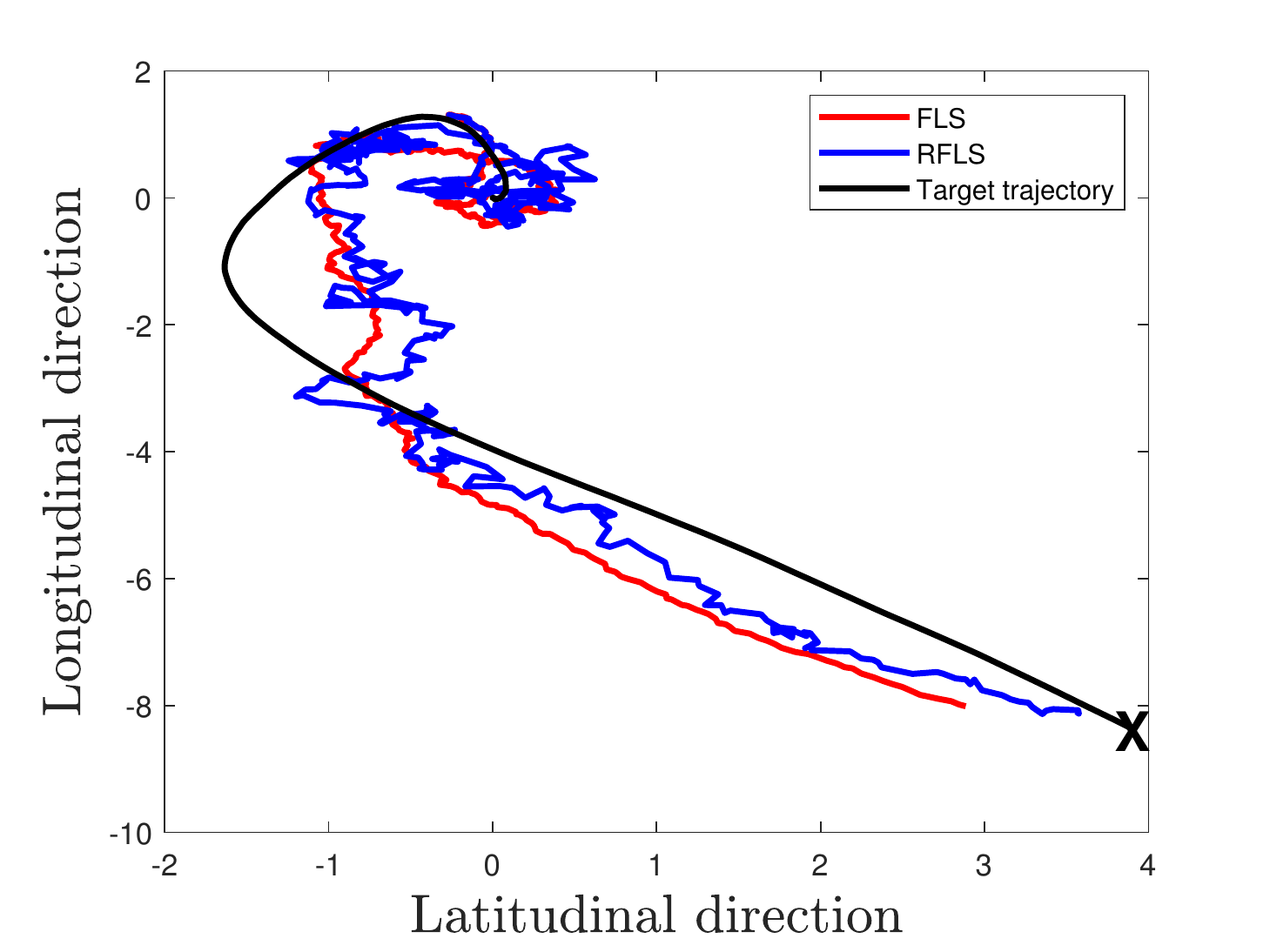}\label{pos_2}}\ \
\vspace{-0.05in}
\caption{ \footnotesize Trajectories generated by the least favorable models in the different ambiguity sets (black line) and the estimated trajectories with RFLS  (blue line) and FLS (red line). The cross denotes the endpoint of the target trajectory.}     \label{pos}
\end{figure}

%

In what follows,  we consider a sample $ Y_N=\{ y_0, y_1\ldots  y_N\}$ with $N=500$ generated by the least favorable model (\ref{LF_x}) in the first ambiguity set with $c_t=10^{-3}$ and  $\xi_0=0$.
\begin{table}
\begin{center}
\setlength{\tabcolsep}{1.8mm}{
\begin{tabular}{ |c|c|c|c|c| }

 \hline
 & RFLS   & FLS & RFLS  & FLS  \\
 & {\scriptsize $c_t=10^{-3}$}&  {\scriptsize$c_t=10^{-3}$}   &  {\scriptsize $c_t=5\cdot 10^{-3}$}   &  {\scriptsize$c_t=5\cdot 10^{-3}$}   \\
 \hline
 {\small $\operatorname{RMSE}_{lat} $} & 0.2918 & 0.4093 & 0.4962 & 0.6112 \\
{\small $\operatorname{RMSE}_{lon}$} & 0.2891 & 0.3280 & 0.4804 & 0.6115\\
 \hline
\end{tabular}}
\end{center}
\caption{RMSE along the latitudinal direction and longitudinal direction.}
\label{tabella}\end{table}
Fig. \ref{pos_1} shows the target trajectory (black line) and the ones estimated with  RFLS (blue line) and  FLS (red line). We also perform the same experiment using the least favorable model in the second ambiguity set with $c_t=5 \times 10^{-3}$.  Moreover, we quantify the performance of the smoothers  in these two experiments through the root mean-square error (RMSE) along the latitudinal direction and longitudinal direction
$$
\operatorname{RMSE}_{lat} = \sqrt{\frac{1}{N} \sum_{t=1}^{N}\left(p_{t}^{lat}-\hat{p}_{t}^{lat}\right)^{2}},
$$
$$
\operatorname{RMSE}_{lon} = \sqrt{\frac{1}{N} \sum_{t=1}^{N}\left(p_{t}^{lon}-\hat{p}_{t}^{lon}\right)^{2} }
$$  whose values are displayed in Table \ref{tabella}.
As expected, RFLS better reduces the influence of the modeling error on the estimation accuracy than others.

\section{Robust parameter estimation} \label{sec_6}
Consider the following state space model
\begin{equation}\label{lf_em}\begin{aligned}
\tilde {\mathcal{M}}(\alpha) : x_{t+1} &={A}(\alpha) x_{t}+{B}(\alpha) \tilde \upsilon_t \\
y_{t} &={C}(\alpha) x_{t}+{{D}(\alpha)} \tilde \upsilon_t
\end{aligned}\end{equation}
where $x_0 \sim \mathcal N(\hat x_0, V_0)$, $\tilde \upsilon_t \sim \mathcal N (0, \tilde R_t)$, i.e. $\tilde v_t$ is a nonstationary process, and the matrices ${A}(\alpha)$, ${B}(\alpha)$, ${C}(\alpha)$ and ${D}(\alpha)$ are parameterized by $\alpha$. In many practical applications, $\alpha$ is not known and needs to be estimated from the observed data $Y_N=\{y_{0}, y_1\ldots  y_{N}\}$. In plain words, the latter is a system identification problem where the model class is $\tilde{\mathcal{M}}=\{\tilde {\mathcal{M}}(\alpha), \alpha \in \Theta \}$ and $\Theta$ is the parameter space. A well established paradigm to find $\alpha$ is the  maximum likelihood (ML) principle. However, it is usually difficult to find an explicit expression of the ML function under Model (\ref{lf_em}). Such a problem is typically addressed by using the expectation-maximization (EM) algorithm, see  \cite[Algorithm 12.3]{sarkka2013bayesian}, which computes a lower bound of the ML function through the iterative scheme:
\begin{itemize}
\item Set up an initial guess $\alpha^0$;
\item For $n=0,1,...$:
    \begin{enumerate}
    \item[$\bullet$] \textbf{E-step}: compute $$ \begin{aligned}
    \mathcal{Q}(\alpha, \alpha^n):=\int  \tilde f_{\alpha^n}(X_{N+1}|Y_N) \log \tilde f_{ \alpha}(Y_N, X_{N+1}) d X_{N+1}
    \end{aligned}$$
    \item[$\bullet$] \textbf{M-step}: compute $\alpha^{(n+1)}=\arg \max_{\alpha} \mathcal{Q}(\alpha, \alpha^n)$
    \end{enumerate}
\end{itemize}
where $\tilde f_{\alpha}(Y_{N}, X_{N+1})$ is the joint density of $Y_N=\{y_0,y_1\ldots,y_N\}$ and $X_{N+1}=\{x_0,x_1,\ldots, x_{N+1}\}$ under $\tilde{\mathcal M}(\alpha) $. In most cases, however,  the covariance matrix $\tilde R_t$ of the noise process  $\tilde v_t$ is not known and it is also time-varying. Such a matrix is typically designed empirically. However, this would require to have the possibility to make more experiments on the system, i.e. a requirement that is not always met. Alternatively, we can select a nominal covariance matrix for the noise process using some a priori knowledge. However, this causes the nominal model to be inconsistent with the actual one and thus the reliability of the estimate of $\alpha$ will be compromised.

A possible way to address this model uncertainty is to understand (\ref{lf_em}) as the least favorable model (\ref{LF_x}) where $\tilde \upsilon_t$ is equal to $\tilde S_t\tilde e_t+L_t\epsilon_{t}$. More precisely, assume that we want to estimate $\alpha$ only knowing the nominal state space model $\mathcal M$, i.e.
\begin{equation}\begin{aligned}\label{iden}
\mathcal{M}(\alpha) : x_{t+1} &={A}(\alpha) x_{t}+{B}(\alpha) \upsilon_t \\
y_{t} &={C}(\alpha) x_{t}+{{D}(\alpha)} \upsilon_t
\end{aligned}\end{equation}
where $v_t$ is normalized WGN, in particular it is a stationary process. It is worth noting that the least favorable model in (\ref{lf_em}) does not belong to $\mathcal M$. Notice that, it is not restrictive to assume that the covariance matrix of $v_t$ is equal to the identity. Indeed, in the case its covariance matrix is $R$, then we can always take $\tilde B(\alpha)=B(\alpha) R^{1/2}$ and $\tilde D(\alpha)=D(\alpha) R^{1/2}$. Then, the least favorable model, solution to (\ref{minimax}), is (\ref{lf_em}) with $\tilde v_t=\tilde S_t\tilde e_t+L_t\epsilon_{t}$. Hence, \begin{equation*} \label{T_R}
\tilde R_t=\tilde S_t \left[\begin{array}{cc} 0 & I_n  \end{array}\right] \Pi_t \left[\begin{array}{c} 0 \\ I_n  \end{array}\right] \tilde S_t^T+L_t L_t^T. \end{equation*}
At this point we can use the density, say $\tilde f^{0}_{\alpha^n}(X_{N+1},Y_N)$, of the least favorable model in order to compute the lower bound $\mathcal Q(\alpha,\alpha^n)$. On the other hand, an approximation of the moments of $\tilde f_{\alpha^n}^0$ required in $\mathcal Q(\alpha,\alpha^n)$ can be constructed by the robust fixed-lag smoother of Algorithms \ref{efficent} and \ref{code:least} leading to
\begin{equation*}
\begin{aligned}
\mathcal{Q}&(\alpha, \alpha^n) \propto -\frac{1}{2} \sum_{t=0}^{N-1} \log |2 \pi B(\alpha)  \tilde R_{t} B^T (\alpha) |\\
&-\frac{1}{2} \tr \left\{  \sum_{t=0}^{N-1} \left[\left( B(\alpha)  \tilde R_{t} B^T (\alpha)\right)^{-1}  \left(\Phi_{1,t} \right.\right. \right. \\
&\left. \left. \left. -\Phi_{2,t} A^T(\alpha) -A(\alpha) \Phi_{2,t}^T + A(\alpha) \Phi_{3,t} A^T(\alpha) \right) \right] \right\}\\
&-\frac{1}{2}  \sum_{t=0}^{N-1} \log |2 \pi D(\alpha) \tilde R_{t+1} D^T (\alpha) |\\
&-\frac{1}{2} \tr \left\{ \sum_{t=0}^{N-1} \left[\left(D(\alpha)  \tilde R_{t+1} D^T (\alpha)\right)^{-1}  \left(\Phi_{4,t} \right. \right. \right.\\
&\left. \left. \left.-C(\alpha)\Phi_{5,t}^T-\Phi_{5,t} C^T(\alpha)+C(\alpha)\Phi_{1,t} C^T(\alpha) \right)\right] \right\}
\end{aligned}\end{equation*}
where
\begin{align*}
\Phi_{1,t}&={\bar \Pi^{L,L}}_{t+L-1}+\hat x_{t+1|t+L-1} \hat x^T_{t+1|t+L-1}, \\\Phi_{2,t}&=({\bar\Pi}^{L,L-1}_{t+L-1})^T+\hat x_{t+1|t+L-1} \hat x^T_{t|t+L-1},\\
\Phi_{3,t}&={\bar\Pi}^{L,L}_{t+L-1}+\hat x_{t|t+L-1} \hat x^T_{t|t+L-1},\\\Phi_{4,t}&=y_{t+1} y^T_{t+1},~~\Phi_{5,t}=y_{t+1} \hat x^T_{t+1|t+L-1},
\end{align*}
and ${\bar \Pi^{j,k}}_{t}$ is the $n\times n$ block in position $(j,k)$ of $\tilde \Pi_{t}$. Then, $\hat x_{t|t+L-1}$, $\hat x_{t-1|t+L-1}$ and $y_t$ are given by Algorithm \ref{efficent}. Clearly, such an approximation is legitimate if the lag $L$ is chosen big enough.

\begin{figure}[htb]
    \centering
    \includegraphics[width=0.65\textwidth]{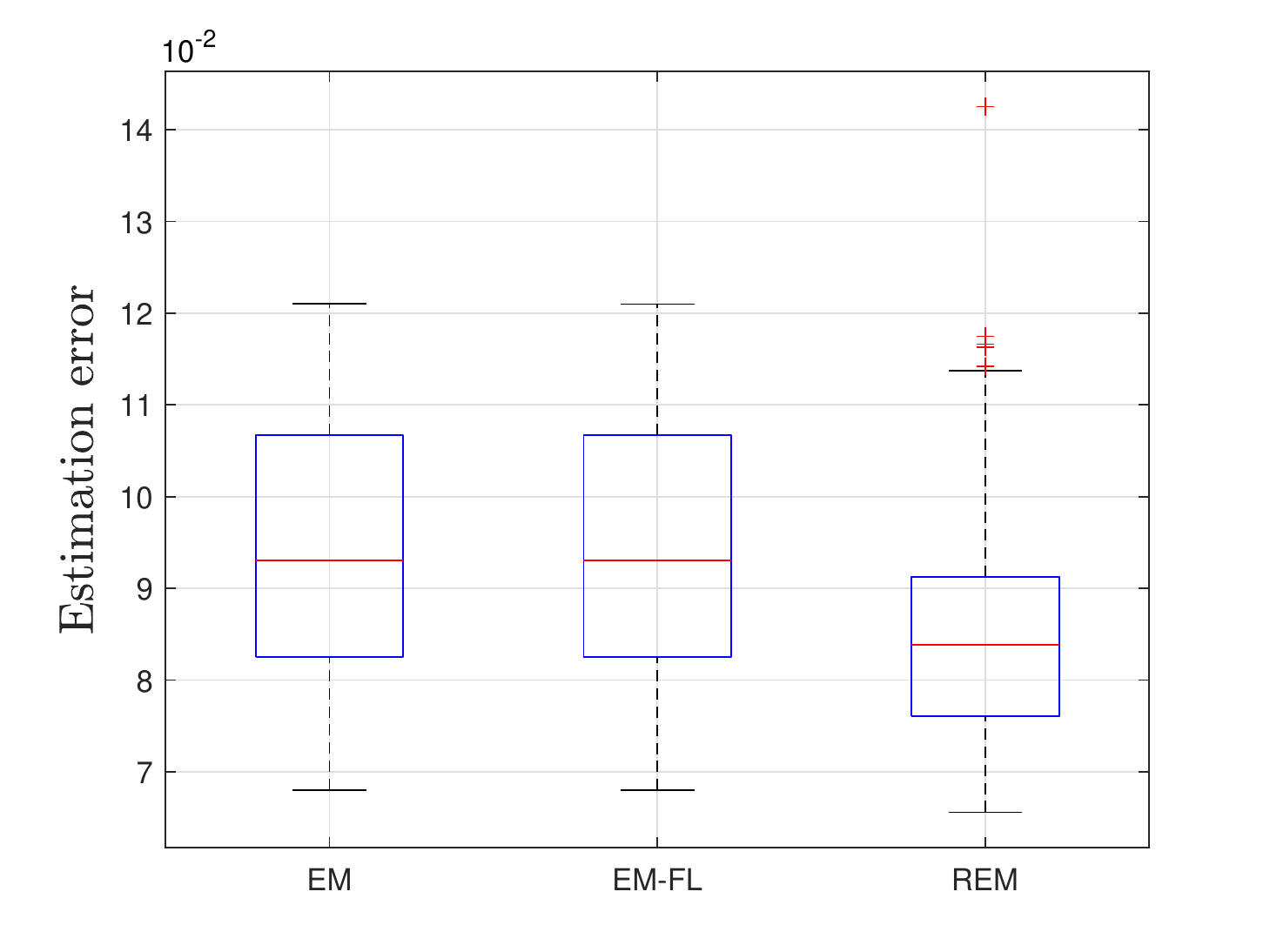}
    \caption{\footnotesize{Estimation error for EM, EM-FL and REM in the  Monte Carlo experiment.}}
    \label{error2}
\end{figure}

Next,  we show a numerical example. We consider the problem to estimate the parameter $\alpha=[\,\alpha_1\; \alpha_2\,]$ using the nominal model class (\ref{iden}) with $$ \begin{aligned}
A(\alpha)&=\left[\begin{array}{cc}
\alpha_1 & 1 \\
0 & \alpha_2
\end{array}\right], ~~~~~ B=\left[\begin{array}{ccc}
0.01 & 0 & 0 \\
0 & 0.45 & 0
\end{array}\right] \\
C&=\left[\begin{array}{ll}
1 & -1
\end{array}\right], ~~~~~~ D=\left[\begin{array}{lll}
0 & 0 & 0.01
\end{array}\right],
\end{aligned}
$$
and the collected data $Y_N=\{y_0,y_1\ldots y_N\}$.
We assume $x_{0} \sim \mathcal{N}\left(0, V_{0}\right)$ with $V_{0}=0.0001 I_{2}$. We assume that the actual model has the same structure of (\ref{iden}) with $\alpha^{\star}=[\,0.1\; 0.9\, ]$, but the actual noise process, say $\tilde v_t$, is not stationary and not known. We solve the aforementioned system identification problem by means of the REM method introduced in Section \ref{sec_6} with $L=50$. Moreover, we compare it with: the standard EM method where the state estimation task is performed by the RTS smoother; the ``fixed-lag'' EM (EM-FL) method where the state estimation task is performed by the standard fixed-lag smoother. To estimate the effectiveness of the REM method we assume that the actual model is the least favorable one belonging to the ambiguity set with $c=2\cdot 10^{-2}$. Moreover, we consider a Monte Carlo experiment with $100$ trials. More precisely, in each trial, we generate the data set $Y_N$ according to (\ref{LF_x}) with $\alpha=\alpha^\star$ and $N=1000$. The initial parameter estimates $\alpha_1^0$ and $\alpha_2^0$, with $\alpha^0:=[\,\alpha_1^0\;\alpha_2^0\,]$, are drawn from a uniform distribution with interval $[0.4,0.9]$ and $[0.07,0.13]$, respectively. Then, the termination condition is $ \left\|\alpha^{n+1}-\alpha^n\right\|\leq \epsilon $ where $\epsilon = 10^{-3}$.

Then, in order to compare the performance of these algorithms, we consider the estimation error $\|\hat\alpha-\alpha^\star\|$ where $\hat\alpha$ is the parameter estimate obtained at the last stage by EM/EM-FL and REM. Fig. \ref{error2} shows the boxplot of the estimation error for EM, EM-FL and REM. We see that REM outperforms EM and EM-FL. It is worth noting that EM and EM-FL perform in the same way. This means that the value of the lag $L$ has been chosen large enough in EM-FL and thus the fixed-lag smoother represents a good approximation of the RTS smoother.

\begin{figure}[t]
\vspace{-0.1in}
\centering
\subfigure[First state component]{\includegraphics[width=2.5in]{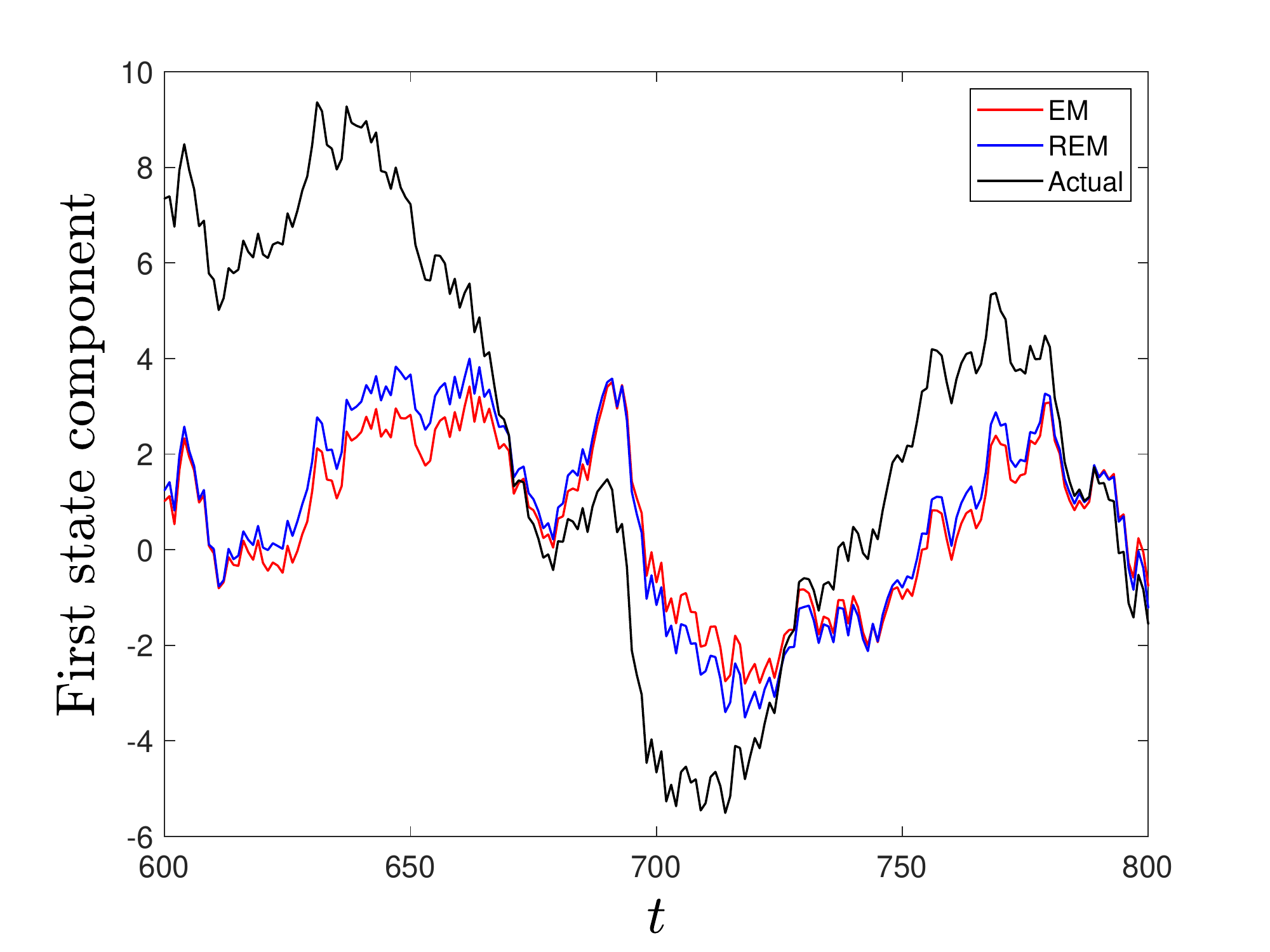}\label{X_1}}\ \
\hspace{3pt}
\subfigure[Second state component]{\includegraphics[width=2.5in]{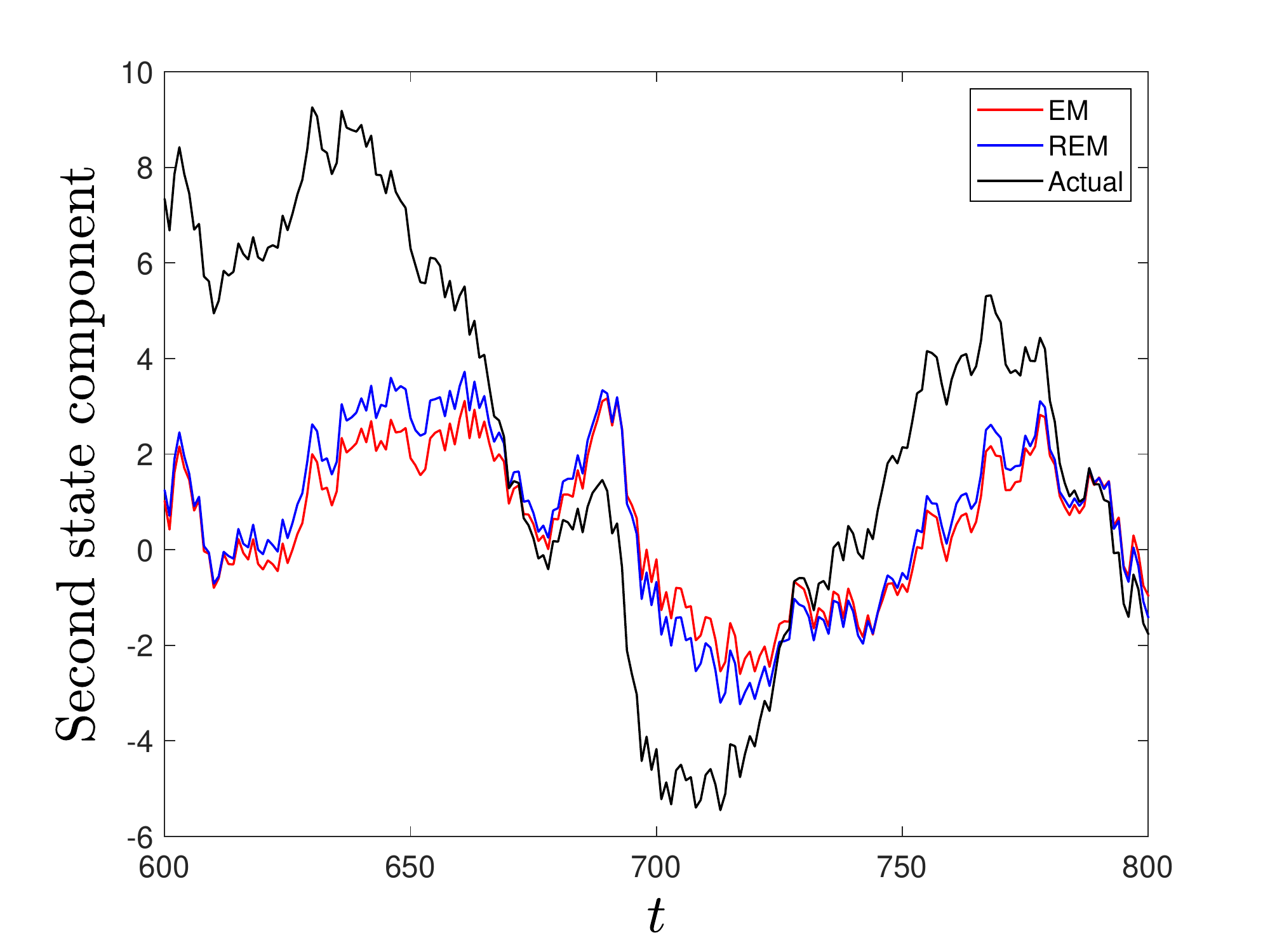}\label{X_2}}\ \
\vspace{-0.05in}
\caption{ \footnotesize Comparison among the actual state components and the estimated ones obtained by EM and REM over the interval $[600, 800]$.}     \label{X}      \vspace{-0.1in}
\end{figure}

Finally, Fig. \ref{X} shows the estimated state trajectory in the last stage by EM and REM,  with $t\in[600, 800]$, in a trial of the Monte Carlo experiment. As we can see, the one obtained with REM, and thus using RFLS, is slightly better than the one with EM, and thus using RTS. Although this advantage is not prominent, it made a dramatic improvement in the performance of the parameter estimator.


\section{Conclusion} \label{sec_8}
In this paper, we have proposed a robust fixed-lag smoother in the case that the actual model is different from the nominal one. More precisely, this paradigm solves a minimax game with two players:  one selects the least favorable model in a prescribed ambiguity set, the other designs the optimal estimator based on the least favorable model. We also proposed an efficient implementation of the robust fixed-lag smoother in order to reduce the computational burden and avoid numerical instabilities. Then, we have characterized the least favorable model for the robust smoother over a finite time horizon. Finally, we have presented some numerical examples showing the effectiveness of the proposed robust fixed-lag smoother.

%


\end{document}